\documentclass[12pt, reqno]{amsart}
\usepackage{ifthen,amsfonts,amsmath,amssymb,epic,eepic}
\usepackage{epsfig,color}

\makeatother

\theoremstyle{definition}

\def\fnum{equation}
\newtheorem{Thm}[\fnum]{Theorem}
\newtheorem{Cor}[\fnum]{Corollary}

\newtheorem{Lem}[\fnum]{Lemma}

\newtheorem{Rem}[\fnum]{Remark}
\newtheorem{Pro}[\fnum]{Proposition}

\numberwithin{equation}{section}

\newcommand{\K}{{\text{K}}}
\newcommand{\nn}{{\bf{n}}}

\newcommand{\dist}{{\text {dist}}}

\newcommand{\an}{{\mathcal{T}}}

\newcommand{\Hess}{{\text {Hess}}}
\def\ZZ{{\bold Z}}
\def\RR{{\bold R}}
\def\SS{{\bold S}}

\def\CC{{\bold C }}
\newcommand{\dd}{\text{r}}

\newcommand{\e}{{\text {e}}}
\newcommand{\Area}{{\text {Area}}}

\newcommand{\tpi}{{\tilde{\Pi}}}

\newcommand{\cC}{{\mathcal{C}}}

\newcommand{\cB}{{\mathcal{B}}}
\newcommand{\Genus}{{\text{gen}}}

\newcommand{\cP}{{\mathcal{P}}}

\newcommand{\cT}{{\mathcal{T}}}

\newcommand{\eqr}[1]{(\ref{#1})}

\begin{document}

\title[Planar domains]
{The space of embedded minimal surfaces of fixed genus in a
$3$-manifold III; Planar domains}

\author{Tobias H. Colding}%
\address{Courant Institute of Mathematical Sciences and Princeton University\\
251 Mercer Street\\ New York, NY 10012 and Fine Hall, Washington Rd.,
Princeton, NJ 08544-1000}
\author{William P. Minicozzi II}%
\address{Department of Mathematics\\
Johns Hopkins University\\
3400 N. Charles St.\\
Baltimore, MD 21218}
\thanks{The first author was partially supported by NSF Grant DMS 9803253
and an Alfred P. Sloan Research Fellowship
and the second author by NSF Grant DMS 9803144
and an Alfred P. Sloan Research Fellowship.}


\email{colding@cims.nyu.edu and minicozz@math.jhu.edu}

\maketitle


\section{Introduction}   \label{s:s0}

This paper is the third in a series  where we describe the space
of all embedded minimal surfaces of fixed genus in a fixed (but
arbitrary) closed $3$-manifold.  In \cite{CM3}--\cite{CM5} we
describe the case where the surfaces are topologically  disks on
any fixed small scale (in fact, Corollary \ref{c:remsingG1} below
is used in \cite{CM5}).  To describe general planar domains (in
\cite{CM6}) we need in addition to the results of
\cite{CM3}--\cite{CM5} a key estimate for embedded stable annuli
which is the main result of this paper (see Theorem
\ref{t:remsing} below).  This estimate asserts that such an
annulus is a graph away from its boundary if it has only one
interior boundary component and if this component lies in a small
(extrinsic) ball.

Planar domains arise when one studies convergence of embedded
minimal surfaces of a fixed genus in a fixed $3$-manifold. This
is due to the next theorem which loosely speaking asserts that any
sequence of embedded minimal surfaces of fixed genus has a
subsequence which are uniformly planar domains away from finitely
many points. (In fact, this describes only ``1).'' and ``2).'' of Theorem
\ref{t:uniplanar}. Case ``3).'' is self explanatory and ``4).'' very
roughly corresponds to whether the surface locally ``looks like''
the genus one helicoid, cf. \cite{HoKrWe}, or has ``more than one end.'')

Given a surface $\Sigma$ with boundary $\partial \Sigma$,
the {\it genus} of $\Sigma$ ($\Genus (\Sigma)$) is the
genus of the closed surface $\hat{\Sigma}$ obtained by adding a
disk to each boundary circle.  The genus of a union of disjoint
surfaces is  the sum of the genuses. Therefore, a
surface with boundary has nonnegative genus; the genus is zero if
and only if it is a planar domain. For example, the disk and the
annulus are both genus zero; on the other hand, a closed surface
of genus $g$ with $k$ disks removed has genus $g$.

\begin{figure}[htbp]
    \setlength{\captionindent}{20pt}
    \begin{minipage}[t]{0.5\textwidth}
    \centering\input{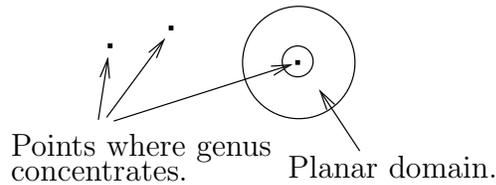}
    \caption{1) and 2) of Theorem  \ref{t:uniplanar}:\newline
    Any sequence of genus $g$ surfaces has a subsequence
    for which the genus
    concentrates at at most $g$ points.\newline
    Away from these points, the surfaces are locally planar
    domains.}
    \label{f:uniplanar}   \end{minipage}
\end{figure}

In the next theorem, $M^3$ will be a closed $3$-manifold and
$\Sigma_i^2$ a sequence of closed embedded oriented
minimal surfaces in $M$
with fixed genus $g$.

\begin{Thm}  \label{t:uniplanar}
See fig. \ref{f:uniplanar}.
There exist $x_1,\dots,x_{m}\in M$ with $m\leq g$ and a subsequence
$\Sigma_j$ so:
\\
\noindent
1). For  $x\in M\setminus \{ x_1,\dots,x_{m} \}$,
there are  $j_x , r_x > 0$ so  $\Genus (B_{r_x}(x)\cap \Sigma_j) =
0$ for  $j > j_x$.
\\
\noindent
2). For each $x_k$, there are $\ell_k , r_k
> 0, r_k
> r_{k,j} \to 0$ so
 for  all $j$
there are components $\{\Sigma^{\ell}_{k,j} \}_{\ell \leq \ell_k}$
 of $B_{r_k} (x_k) \cap \Sigma_j$ with
 $\Genus(B_{r_k}(x_k) \cap \Sigma_j) = \sum_{\ell \leq \ell_k} \Genus(
\Sigma^{\ell}_{k,j} )\leq g$ and
$\Genus(\Sigma_{k,j}^{\ell}) =
\Genus(B_{r_{k,j}}(x_k) \cap \Sigma_{k,j}^{\ell})$ for $\ell \leq \ell_k$.
\\
\noindent
3).  For   every $k , \ell , j$, there is only one component
$\tilde \Sigma_{k,j}^{\ell}$
 of $B_{r_{k,j}}(x_k) \cap \Sigma_{k,j}^{\ell}$ with genus $> 0$.
\\
\noindent
4). For each $k,\ell$, either $\partial \Sigma_{k,j}^{\ell}$ is connected or a
component of $\partial \tilde \Sigma_{k,j}^{\ell}$ separates two
components of $\partial  \Sigma_{k,j}^{\ell}$.
\end{Thm}

\begin{figure}[htbp]
    \setlength{\captionindent}{20pt}
    \begin{minipage}[t]{0.5\textwidth}
    \centering\input{pl2a.pstex_t}
    \caption{The catenoid given by revolving $x_1= \cosh x_3$
around the $x_3$-axis.}
    \label{f:cat}   \end{minipage}\begin{minipage}[t]{0.5\textwidth}
    \centering\input{pl2.pstex_t}
    \caption{The Riemann examples: Parallel planes connected by necks.}
    \label{f:riemann}   \end{minipage}
\end{figure}

To explain why the next two theorems are crucial for
what we call ``the pair of pants decomposition'' of embedded minimal
planar domains, recall the following prime examples of such domains:
Minimal graphs (over disks),
a helicoid, a
catenoid or one of the Riemann examples. (Note that the first two are
topologically disks and the others are disks with one or more subdisks
removed.)  Let us describe the non simply connected
examples in a little more detail.  The catenoid (see fig. \ref{f:cat})
is the (topological) annulus
\begin{equation}
        (\cosh s\, \cos t,\cosh s\, \sin t,s)
\end{equation}
where $s,\,t\in\RR$.  To describe the Riemann examples, think of a
catenoid as roughly being obtained by connecting two parallel
planes by a neck.  Loosely speaking (see fig. \ref{f:riemann}),
the Riemann examples are
given by connecting (infinitely many) parallel planes by necks;
each adjacent pair of planes is connected by exactly one neck.  In
addition, all of the necks are lined up along an axis and the
separation between each pair of adjacent ends is constant (in fact
the surfaces are periodic). Locally, one can imagine connecting
$\ell -1$ planes by $\ell-2$ necks and add half of a catenoid to
each of the two outermost planes, possibly with some restriction
on how the necks line up and on the separation of the planes; see
\cite{FrMe}, \cite{Ka}, \cite{LoRo}.

\begin{figure}[htbp]
    \setlength{\captionindent}{20pt}
    \begin{minipage}[t]{0.5\textwidth}
    \centering\input{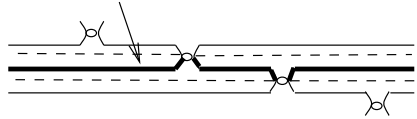}
    \caption{Decomposing the Riemann examples into ``pairs of pants'' by
    cutting along small curves; these curves bound minimal graphical
    annuli separating
    the ends.}
    \label{f:cutriemann}   \end{minipage}
\end{figure}

To illustrate how Theorem \ref{t:remsing} below will be used in
\cite{CM6} where we give the actual ``pair of pants
decomposition'' observe that the catenoid can be decomposed into
two minimal annuli each with one exterior convex boundary and one
interior boundary which is a short simple closed geodesic. (See
also \cite{CM9} for the ``pair of pants decomposition'' in the
special case of annuli.)  In the case of the Riemann examples
(see fig. \ref{f:cutriemann}),
there will be a number of ``pair of pants'', that is, topological
disks with two subdisks removed. Metrically these ``pair of
pants'' have one convex outer boundary and two interior boundaries
each of which is a simple closed geodesic. Note also that this
decomposition can be made by putting in minimal graphical annuli
in the complement of the domains (in $\RR^3$) which separate each
of the pieces; cf. Corollary \ref{c:remsing} below.  Moreover,
after the decomposition is made then every intersection of one of
the ``pair of pants'' with an extrinsic ball away from the
interior boundaries is simply connected and hence the results of
\cite{CM3}--\cite{CM5} apply there.

\begin{figure}[htbp]
    \setlength{\captionindent}{20pt}
    \begin{minipage}[t]{0.5\textwidth}
    \centering\input{pl4.pstex_t}
    \caption{Theorem \ref{t:remsing}: Embedded stable annuli with small interior boundary
    are graphical away from their boundary.}
    \label{f:remsing}   \end{minipage}
\end{figure}

The next theorem is a kind of effective removable singularity
theorem for embedded stable minimal surfaces with small interior
boundaries.   It asserts that embedded stable minimal surfaces
with small interior boundaries are graphical away from the
boundary.  Here small means contained in a small ball in $\RR^3$
(and not that the interior boundary has small length).  This
distinction is important; in particular if one had a bound for the
area of a tubular neighborhood of the interior boundary, then
Theorem \ref{t:remsing} would  follow easily;
see Corollary \ref{c:asmall}
and cf. \cite{Fi}.

\begin{Thm}  \label{t:remsing}
See fig. \ref{f:remsing}.
Given $\tau > 0$, there exists $C_1>1$, so if $\Gamma\subset
B_{R}\subset \RR^3$ is an embedded stable minimal annulus
with $\partial \Gamma \subset \partial B_{R} \cup
B_{r_0/4}$ (for $C_1^2 \, r_0 < R$) and $B_{r_0}\cap \partial \Gamma$
is connected, then each component of
$B_{R / C_1} \cap \Gamma \setminus
B_{C_1 \, r_0}$ is a graph with gradient $\leq \tau$.
\end{Thm}

Many of the results of this paper will involve either graphs or
multi-valued graphs.  Graphs will always be assumed to be
single-valued over a domain in the plane (as is the case in
Theorem \ref{t:remsing}).

 For simplicity, $\Gamma$ in Theorem \ref{t:remsing} is assumed
to be an annulus;
see \cite{CM6} for the slight additional arguments needed when
$\Gamma$ is a general planar domain.  (Once we see in \cite{CM6} that
Theorem \ref{t:remsing} holds for general planar domains then the
corresponding generalization of Corollary \ref{c:remsing} follows.)

\begin{figure}[htbp]
    \setlength{\captionindent}{20pt}
    \begin{minipage}[t]{0.5\textwidth}
    \centering\input{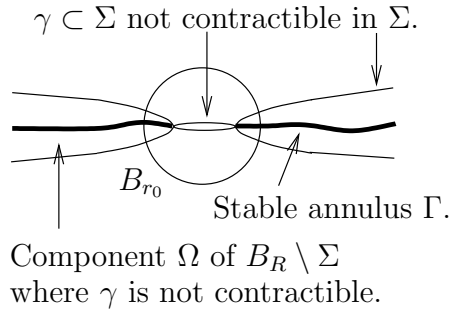}
    \caption{Corollary \ref{c:remsing}: Solving a Plateau problem gives
    a stable graphical annulus separating the boundary components
    of an embedded minimal annulus.}
    \label{f:barrier}   \end{minipage}
\end{figure}

 Combining Theorem \ref{t:remsing} with the solution of a Plateau
problem of Meeks-Yau (proven initially for convex domains in
theorem 5 of \cite{MeYa1} and extended to mean convex domains in
\cite{MeYa2}), we get (the result of Meeks-Yau gives the existence
of $\Gamma$ below):

\begin{Cor} \label{c:remsing}
See fig. \ref{f:barrier}.
Given $\tau > 0$, there exists $C_1>1$, so: Let $\Sigma\subset
B_{R}\subset \RR^3$, $\partial \Sigma\subset \partial B_R$ be an
embedded minimal surface with $\Genus (\Sigma)=\Genus (B_{r_1}\cap
\Sigma)$ and let $\Omega$ be a component of $B_R \setminus
\Sigma$. If $\gamma\subset B_{r_0} \cap \Sigma\setminus B_{r_1}$
is noncontractible and homologous in $\Sigma \setminus B_{r_1}$ to
a component of $\partial\Sigma$ and $r_0>r_1$, then a component
$\hat{\Sigma}$ of $\Sigma \setminus \gamma$ is an annulus and
there is a stable embedded minimal annulus $\Gamma \subset \Omega$
with $\partial \Gamma = \partial \hat{\Sigma}$. Moreover,  each
component of $(B_{R/C_1}\setminus B_{C_1 \, r_0})\cap \Gamma$ is a
graph with gradient $\leq \tau$.
\end{Cor}

Stability of $\Gamma$ in Theorem \ref{t:remsing} is used in two
ways: To get a pointwise curvature bound on  $\Gamma$ and to show
that certain sectors have small curvature.  In section $2$
of \cite{CM4}, we showed that a pointwise curvature bound allows
us to decompose an embedded minimal surface into a set of bounded
area and a collection of (almost stable) sectors with small curvature.  Using this, the proof of Theorem
\ref{t:remsing} will also give
(if $0\in \Sigma$, then
$\Sigma_{0,t}$ denotes the component of $B_t \cap
\Sigma$ containing $0$):

\begin{Thm}  \label{t:remsingG}
Given $C$, there exist $C_2 , C_3 > 1$, so: Let $0 \in
\Sigma\subset B_{R}\subset \RR^3$ be an embedded minimal surface
with connected $\partial \Sigma \subset \partial B_{R}$. If
$\Genus (\Sigma_{0,r_0}) = \Genus ( \Sigma )$, $r_0 \leq R/ C_2$,
and
\begin{equation}    \label{e:atcbG}
\sup_{\Sigma\setminus B_{r_0}} |x|^2 \, |A|^2(x) \leq C \, ,
\end{equation}
then $\Area (\Sigma_{0,r_0} ) \leq  C_3 \, r_0^2$.
\end{Thm}

In \cite{CM5} a strengthening of Theorem \ref{t:remsingG} (this
strengthening is Theorem \ref{t:remsingG1} below) will be used to
show that for limits of a degenerating sequence of embedded
minimal disks points where the curvatures blow up are not
isolated.  This will eventually give (theorem $0.1$ of \cite{CM5})
that for a subsequence such points form a Lipschitz curve which is
infinite in two directions and transversal to the limit leaves;
cf. with the example given by a sequence of rescaled helicoids.

To describe a neighborhood of each of the finitely many points,
coming from Theorem \ref{t:uniplanar},
where the genus concentrates (specifically to describe when
there is one component $\tilde \Sigma_{k,j}^{\ell}$ of genus $>0$ in ``3).''
of Theorem \ref{t:uniplanar}), we will need in \cite{CM6}:

\begin{Cor} \label{c:d}
Given $C,g$, there exist $C_4, C_5$ so:
Let $0 \in \Sigma\subset B_{R} \subset \RR^3$  be an embedded minimal surface
with connected $\partial \Sigma \subset \partial B_{R}$, $r_0 < R/C_4$,
and $\Genus(\Sigma_{0,r_0}  ) =\Genus(\Sigma) \leq g$.
If
\begin{equation}    \label{e:atcb}
\sup_{\Sigma \setminus B_{r_0}} |x|^2 \, |A|^2(x) \leq C \, ,
\end{equation}
then $\Sigma$ is a disk and $\Sigma_{0,R/C_5}$  is a graph with
gradient $\leq 1$.
\end{Cor}

This corollary follows directly from Theorem \ref{t:remsingG} and
theorem 1.22 of \cite{CM4}. Namely, note first that for $r_0 \leq s \leq R$,
  it follows from the
maximum principle (since $\Sigma$ is minimal) and Corollary
\ref{c:gen2} that $\partial \Sigma_{0,s}$ is connected and $\Sigma
\setminus \Sigma_{0,s}$ is an annulus.
Second, Theorem \ref{t:remsingG}
bounds  $\Area (\Sigma_{0,R/C_2})$ and theorem 1.22 of \cite{CM4}
then gives the corollary.

Theorems \ref{t:remsing}, \ref{t:remsingG} and Corollary \ref{c:d}
are local and are for simplicity
stated and proven only in $\RR^3$ although
they can with only very minor changes easily be seen to hold for
minimal planar domains in a sufficiently small ball in any
given fixed Riemannian $3$-manifold.

\vskip2mm
Throughout $\Sigma$, $\Gamma\subset M^3$ will denote
complete minimal surfaces possibly with boundary, sectional
curvatures $\K_{\Sigma}$, $\K_{\Gamma}$, and second fundamental
forms $A_{\Sigma}$, $A_{\Gamma}$.  $\Gamma$ will be assumed to be
stable and have trivial normal bundle. Given $x \in M$, $B_s (x)$
will be the extrinsic geodesic ball with radius $s$ and center
$x$. Likewise, if $x \in \Sigma$, then $\cB_s(x)$ is the intrinsic
ball in $\Sigma$. Given $S\subset \Sigma$ and $t>0$, let
$\an_{t}(S,\Sigma) \subset \Sigma$ be the intrinsic tubular
neighborhood of $S$ in $\Sigma$ with radius $t$ and set
$\an_{s,t}(S,\Sigma) = \an_{t}(S,\Sigma) \setminus \an_s
(S,\Sigma)$.  Unless explicitly stated otherwise, all geodesics
will be parametrized by arclength.

We will often consider the intersections of various curves and surfaces
with extrinsic balls.  We will always assume that
these intersections are transverse since this can anyway be achieved
by an arbitrarily small perturbation of the radius.

\setcounter{part}{0}
\numberwithin{section}{part} 
\renewcommand{\rm}{\normalshape} 
\renewcommand{\thepart}{\Roman{part}}
\setcounter{section}{1}

\part{Topological decomposition of surfaces}  \label{p:p0}

\begin{figure}[htbp]
    \setlength{\captionindent}{20pt}
    \begin{minipage}[t]{0.5\textwidth}
    \centering\input{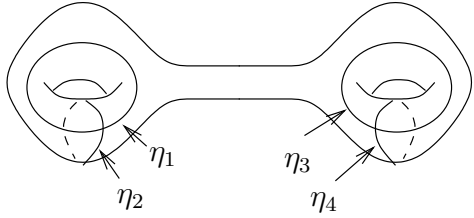}
    \caption{Lemma \ref{l:gen1}: A basis for homology on a
    surface of genus $g$.}
    \label{f:gen1}   \end{minipage}
\end{figure}

In this part we will first collect some simple facts and results
about planar domains
and domains that are planar outside a small ball.  These results will then
be used to show Theorem \ref{t:uniplanar}.  First we have:

\begin{Lem}    \label{l:gen1}
See fig. \ref{f:gen1}.
Let $\Sigma$ be a closed oriented surface (i.e., $\partial
\Sigma = \emptyset$) with genus $g$.  There are transverse
 simple closed curves
$\eta_1 , \dots , \eta_{2g} \subset \Sigma$ so that for $i <
j$
 \begin{equation}   \label{e:scc2}
 \# \{ p \, | \, p \in \eta_i  \cap \eta_j \} =  \delta_{i+g,j} \, .
\end{equation}
Furthermore, for any such $\{ \eta_i \}$,  if $\eta \subset
\Sigma \setminus \cup_i \eta_i$ is a  closed curve, then
$\eta$ divides $\Sigma$.
\end{Lem}

Recall that if $\partial \Sigma \ne \emptyset$, then $\hat \Sigma$ is
the surface obtained by replacing each circle in $\partial \Sigma$
with a disk.
Note that a closed curve $\eta \subset \Sigma$ divides $\Sigma$ if
and only if $\eta$ is homologically trivial in $\hat{\Sigma}$.

\begin{Cor}    \label{c:gen2}
If $\Sigma_1 \subset \Sigma$  and $\Genus ( \Sigma_1) = \Genus
(\Sigma)$, then each simple closed curve $\eta \subset \Sigma
\setminus \Sigma_1$ divides $\Sigma$.
\end{Cor}

\begin{proof}
Since $\Sigma_1$ has genus $g= \Genus (\Sigma)$, Lemma
\ref{l:gen1} gives transverse simple closed  curves $\eta_1 ,
\dots , \eta_{2g} \subset  \Sigma_1$ satisfying \eqr{e:scc2}.
However, since $\eta$ does not intersect any of the $\eta_i$'s,
Lemma \ref{l:gen1} implies that $\eta$ divides $\Sigma$.
\end{proof}

\begin{Cor} \label{c:thegen}
 If  $\Sigma$ has a
decomposition  $\Sigma =\cup_{\beta=1}^{\ell}\Sigma_{\beta}$ where
the union is taken over the boundaries and each $\Sigma_{\beta}$
is a surface  with boundary consisting of a number of disjoint
circles, then
\begin{equation}    \label{e:lsc}
\sum_{\beta=1}^{\ell} \Genus (\Sigma_{\beta}) \leq \Genus (\Sigma)
\, .
\end{equation}
\end{Cor}

\begin{proof}
Set $g_{\beta} = \Genus (\Sigma_{\beta})$.  Lemma \ref{l:gen1},
gives transverse simple closed curves
  $ \eta^{\beta}_1 , \dots ,  \eta^{\beta}_{2g_{\beta}} \subset
   \Sigma_{\beta}$
satisfying \eqr{e:scc2}. Since $\Sigma_{\beta_1} \cap
\Sigma_{\beta_2} = \emptyset$ for $\beta_1 \ne \beta_2$,   this
implies that the rank of the intersection form on the first
homology (mod $2$) of $\hat{\Sigma}$ is $\geq
2\sum_{\beta=1}^{\ell}g_{\beta}$.  In particular, we get
\eqr{e:lsc}.
\end{proof}

In the next lemma,  $M^3$ will be a closed $3$-manifold and
$\Sigma_i^2$ a sequence of closed embedded oriented
minimal surfaces in $M$
with fixed genus $g$.

\begin{Lem}  \label{l:uniplanar}
There exist $x_1,\dots,x_{m}\in M$ with $m \leq g$ and a subsequence $\Sigma_j$
so: For  $x\in M\setminus \{ x_1,\dots,x_{m} \}$,
there exist  $j_x , r_x > 0$ so  $\Genus (B_{r_x}(x)\cap \Sigma_j) =
0$ for $j > j_x$.   For each
$x_k$, there exist $R_k , g_k
> 0, R_k  > R_{k,j} \to 0$ so
 $\Genus(B_{R_k}(x_k) \cap \Sigma_j) = g_k = \Genus(
 B_{R_{k,j}}(x_k) \cap \Sigma_j )$ for  all $j$ and $\sum_{k=1}^m g_k \leq g$.
\end{Lem}

\begin{proof}
Suppose that for some $x_1\in M$ and any $R_1
> 0$ we have $\Genus (B_{R_1}(x_1) \cap \Sigma_i) = g_{1,i} > 0$ for
infinitely many $i$'s.  By Corollary \ref{c:thegen},  $g_{1,i} \leq g$
and so there is a subsequence $\Sigma_j$ and
a sequence $R_{1,j} \to 0$ so
that for all $j$
\begin{equation}    \label{e:r1x1}
    \Genus (B_{R_{1,j}}(x_1) \cap \Sigma_j) = g_1 > 0 \, .
\end{equation}
By repeating this construction, we can suppose  that there are
disjoint points $x_1 , \dots , x_m \in M$ and $R_{k,j} > 0$ so
that for any $k$ we have $R_{k,j} \to 0$ and
 $\Genus (B_{R_{k,j}}(x_k) \cap \Sigma_j) = g_k > 0$.

However,  Corollary \ref{c:thegen} implies that for $j$
sufficiently large
\begin{equation}
0\leq  \Genus ( \Sigma_j \setminus \cup_k B_{R_{k,j}}(x_k) ) \leq
\Genus (\Sigma_j) - \sum_{k=1}^m \Genus (B_{R_{k,j}}(x_k) \cap
\Sigma_j) \leq g - \sum_{k=1}^m g_k \, .
\end{equation}
In particular, $\sum_{k=1}^m g_k \leq g$ and we can therefore
assume that $\sum_{k=1}^m g_k$ is maximal.
This has two
consequences.  First,  given $x \in M \setminus \{ x_1 , \dots ,
x_m \}$, there exist $r_x
> 0$ and $j_x$ so that $\Genus (B_{r_x}(x) \cap \Sigma_j) = 0$ for
$j
> j_x$.  Second, for each $x_k$, there exist $R_k > 0$ and $j_k$ so
$\Genus(B_{R_k}(x_k) \cap \Sigma_j) = g_k$  for  $j > j_k$.
The lemma now follows easily.
\end{proof}

By Corollary \ref{c:thegen}, each $R_k , R_{k,j}$ from Lemma
\ref{l:uniplanar} can (after going to a further subsequence) be
replaced by any $R_k' , R_{k,j}'$ with
$R_k' \leq R_k$ and $R_{k,j}' \geq R_{k,j}$.
Similarly, each $r_x$ can be replaced by any $r_x' \leq r_x$. This
will be used freely in the proof of Theorem \ref{t:uniplanar}
below.

\begin{proof}
(of Theorem \ref{t:uniplanar}). Let $x_k , g_k , R_k ,
R_{k,j}$ and $r_x$ be from Lemma \ref{l:uniplanar}. We can assume
that each $R_k > 0$ is sufficiently small so that $B_{R_k} (x_k)$
is essentially Euclidean  (e.g., $R_k < \min \{ i_0 / 4 , \pi /
(4 k^{1/2}) \}$).
``1).'' follows directly from Lemma \ref{l:uniplanar}.

  For each $x_k$, we can assume that there are $\ell_k$ and
  $n_{\ell ,k}$ so:
   $B_{R_k}(x_k) \cap \Sigma_j$ has components
   $\{ \Sigma_{k,j}^{\ell} \}_{1\leq \ell \leq
  \ell_k}$
  with genus $> 0$ and  $B_{R_{k,j}}(x_k) \cap \Sigma_{k,j}^{\ell}$ has
$n_{\ell,k}$ components  with genus $>0$.
We will use repeatedly that, by
``1).'' and Corollary \ref{c:thegen}:
$n_{\ell, k}$ is nonincreasing if  either $R_{k,j}$ increases or
$R_k$ decreases.
  For each $\ell ,k$ with $n_{\ell , k} > 1$, set
\begin{equation}
    \rho_{k,j}^{\ell} = \inf \{ \rho >  R_{k,j} \, | \, \# \{ {\text{components of }}
    B_{\rho}(x_k) \cap \Sigma_{k,j}^{\ell} \}  < n_{\ell , k}  \} \, .
\end{equation}
There are two cases.  If $\liminf_{j\to \infty} \rho_{k,j}^{\ell} = 0$, then
choose a subsequence $\Sigma_j$ with
$\rho_{k,j}^{\ell} \to 0$;  $n_{\ell,k}$ decreases if we replace
$R_{k,j}$ with any $R_{k,j}' > \rho_{k,j}^{\ell}$.
Otherwise,   set
$ 2\, \rho_{k}^{\ell} =  \liminf_{j \to \infty} \rho_{k,j}^{\ell} > 0$
and  choose a
subsequence $\Sigma_j$ so $\rho_{k,j}^{\ell} < \rho_k^{\ell}$;
$\ell_k$ increases if we replace
$R_k$ with any $R_k' \leq \rho_k^{\ell}$.
In either case, $\sum_{\ell , k} (n_{\ell,k} - 1)$ decreases.
Since $\sum_{\ell , k} n_{\ell , k} \leq g$  (by Corollary \ref{c:thegen}),  repeating this
$\leq g$ times gives
   $0 < R_k' \leq R_k , R_{k,j} \leq R_{k,j}' \to 0$
   and a subsequence so
only one component $\tilde \Sigma_{k,j}^{\ell}$
 of $B_{R_{k,j}'}(x_k) \cap \Sigma_{k,j}^{\ell}$ has genus $> 0$
   (i.e., each new $n_{\ell ,k}=1$).
By Corollary \ref{c:thegen} (and ``1).'') and the remarks before the proof,
``1).'', ``2).'', and ``3).'' now hold for
any $r_k \leq R_k'$ and $R_{k,j}' \leq r_{k,j}
\to 0$.

Suppose that for some $k,\ell$ there exists $j_{k,\ell}$ so
$\partial \Sigma_{k,j}^{\ell}$ has at least two components for all
$j > j_{k,\ell}$.   For  $R_{k,j}' \leq t \leq R_k'$,
let $\Sigma_{k,j}^{\ell}(t)$ be the component of
$B_t(x_k) \cap \Sigma$ containing $\tilde{\Sigma}_{k,j}^{\ell}$.
Set
\begin{equation}
    r_{k,j}^{\ell} = \inf \{ t >  R_{k,j} \, | \, \# \{ {\text{components of }}
    \partial \Sigma_{k,j}^{\ell}(t)  \}  > 1 \} \, .
\end{equation}
 There are two cases.
If $\liminf_{j\to \infty} r_{k,j}^{\ell} = 0$, then
choose a subsequence $\Sigma_j$ with
$r_{k,j}^{\ell} \to 0$.   By the
maximum principle (since $\Sigma$ is minimal) and Corollary
\ref{c:gen2}, a component of (the new) $\partial \tilde
\Sigma_{k,j}^{\ell}$  separates two components of $\partial
\Sigma_{k,j}^{\ell}$ for any $r_{k,j} \to 0$ with
$r_{k,j} > r_{k,j}^{\ell}$.
On the other hand, if
$\liminf_{j\to \infty} r_{k,j}^{\ell} = 2 r_k^{\ell} > 0$,
then choose a subsequence so
(the new) $\partial \Sigma_{k,j}^{\ell}$  is connected
for any $r_k \leq r_{k}^{\ell}$.
After repeating this $\leq g$ times (each time either
increasing $R_{k,j}'$ or decreasing $R_k'$), ``4).'' also holds.
\end{proof}

In \cite{CM6} we will need the following
(here, and elsewhere, if $0\in \Sigma\subset \RR^3$,
then $\Sigma_{0,t}$ denotes the component of $B_t\cap \Sigma$ containing $0$):

\begin{Pro}
Let $0\in\Sigma_i\subset B_{S_i}\subset \RR^3$ with $\partial
\Sigma_i\subset \partial B_{S_i}$ be a sequence of embedded
minimal surfaces with genus $\leq g<\infty$ and $S_i\to \infty$.
After going to a subsequence, $\Sigma_j$, and possibly replacing
$S_j$ by $R_j$ and $\Sigma_j$ by $\Sigma_{0,j,R_j}$ where $R_0\leq
R_j\leq S_j$ and $R_j\to \infty$, then $\Genus
(\Sigma_{j,0,R_0})=\Genus (\Sigma_j)\leq g$ and either (a)
or (b) holds:\\
(a) $\partial \Sigma_{j,0,t}$ is connected for all $R_0\leq t\leq R_j$.\\
(b) $\partial \Sigma_{j,0,R_0}$ is disconnected.
\end{Pro}

\begin{proof}
We will first show that there exists $R_0 >0$, a subsequence
$\Sigma_j$, and a sequence $R_j\to \infty$ with $R\leq R_j\leq
S_j$, such that (after replacing $\Sigma_j$ by $\Sigma_{j,0,R_j}$)
$\Genus (\Sigma_{j,0,R_0})=\Genus (\Sigma_j)\leq g$.  Suppose not;
it follows easily from the monotonicity of the genus (i.e.,
Corollary  \ref{c:thegen}) that there exists a subsequence
$\Sigma_j$ and a sequence $G_k\to\infty$ such that for all $k$
there exists a $j_k$ so for $j\geq j_k$
\begin{equation}
g \geq \Genus (\Sigma_{j,0,G_{k+1}})>\Genus (\Sigma_{j,0,G_{k}})
    \, ,
\end{equation}
which is a contradiction.

 For each $j$, let $R_{0,j}$ be the infimum of $R$ with $R_0 \leq R
 \leq R_j$ where $\partial \Sigma_{j,0,R}$ is disconnected;
set $R_{0,j} = R_j$ if
 no such exists.
Either $\liminf R_{0,j} < \infty$,
 in which case, after going to a subsequence and replacing $R_0$ by
$\liminf R_{0,j}$+1, we are in
 (b) by the maximum principle.  Or, if $\liminf R_{0,j} =
 \infty$, then we are in (a) after replacing $R_j$ by $R_{0,j}$.
\end{proof}

\part{Estimates for stable minimal surfaces with small interior boundaries}
\label{p:p1}

In this part we prove Theorem \ref{t:remsing}.

\section{Long stable sectors contain multi-valued graphs}   \label{s:lst}

In \cite{CM3}, \cite{CM4} we gave estimates for stable sectors. A
stable sector in the sense of \cite{CM3}, \cite{CM4} is a stable
subset of a minimal surface given as half of a normal tubular
neighborhood (in the surface) of a strictly convex curve (for
instance, a curve lying in the boundary of an intrinsic ball). In
this section we give similar estimates for half of normal tubular
neighborhoods of curves lying in the intersection of the surface
and the boundary of an extrinsic ball. These domains arise
naturally in our main result and are unfortunately somewhat more
complicated to deal with due to the lack of convexity of the
curves.

In this section, the surfaces $\Sigma$ and $\Gamma$
will be
 planar domains and, hence,  simple closed curves will
divide the surface into two planar (sub)domains.

We will need some notation for
multi-valued graphs.  Let $\cP$ be the universal cover of the
punctured plane $\CC \setminus \{ 0 \}$ with global (polar)
coordinates $(\rho , \theta)$ and set $S_{r,s}^{\theta_1 , \theta_2} =
\{  r \leq  \rho \leq  s  \, ,
\, \theta_1 \leq \theta \leq \theta_2 \}$.
An $N$-valued graph $\Sigma$
of a function $u$ over
the annulus $D_{s} \setminus D_{r}$
(see fig. \ref{f:separationw})
is a (single-valued) graph (of $u$)
over $S_{r,s}^{-N\,\pi,N\,\pi}$
($\Sigma_{r,s}^{\theta_1 , \theta_2}$
will denote the subgraph of $\Sigma$
over $S_{r,s}^{\theta_1 , \theta_2}$).
The separation $w(\rho,\theta)$ between consecutive sheets is
(see fig. \ref{f:separationw})
\begin{equation}        \label{e:thesepw}
w(\rho,\theta)=u(\rho,\theta+2\pi)-u(\rho,\theta)\, .
\end{equation}

\begin{figure}[htbp]
    \setlength{\captionindent}{20pt}
    \begin{minipage}[t]{0.5\textwidth}
    \centering\input{pl8a.pstex_t}
    \caption{The separation $w$
for a multi-valued  graph in \eqr{e:thesepw}.}
    \label{f:separationw}   \end{minipage}\begin{minipage}[t]{0.5\textwidth}
    \centering\input{pl8b.pstex_t}
    \caption{Theorem \ref{t:longstable}:
Embedded stable annuli with small interior boundary
    contain either (1) a graphical annulus or (2) an $N$-valued
 graph away from its boundary.}
    \label{f:longstable}   \end{minipage}
\end{figure}

The main result of the next two sections is
($\Gamma_1 (\partial) $  is the component of
$B_{1} \cap \Gamma$ containing $B_1 \cap
\partial \Gamma$):

\begin{Thm}   \label{t:longstable}
See fig. \ref{f:longstable}.
Given $N, \tau > 0$, there exist $\omega > 1$, $d_0$ so: Let
$\Gamma$ be a stable embedded minimal annulus with $\partial
\Gamma \subset  B_{1/4} \cup \partial B_{R}$,
$B_{1/4} \cap \partial \Gamma$  connected,  and $R > \omega^2$.
Given $z_1 \in \partial B_1 \cap
\partial \Gamma_1 (\partial)$,
(after a rotation of $\RR^3$) either (1) or (2) below holds:\\
(1) Each component of $B_{R/\omega} \cap \Gamma \setminus B_{\omega}$
is a graph with gradient $\leq \tau$.\\
(2)  $\Gamma$ contains a graph $\Gamma_{\omega ,
R/\omega}^{-N\pi,N\pi}$ with gradient $\leq \tau$ and
$\dist_{\Gamma\setminus \Gamma_{1}(\partial)} ( z_1 ,
\Gamma_{\omega , \omega }^{0,0}) < d_0$.
\end{Thm}

Note that if $\Gamma$ is as in Theorem \ref{t:longstable}
and one component of  $B_{R/\omega} \cap \Gamma \setminus B_{\omega}$ contains a graph over
$D_{R/(2\omega)} \setminus
D_{2 \omega}$ with gradient $\leq 1$, then every component of
$B_{R/(C\omega)} \cap \Gamma \setminus
B_{C\omega}$ is a graph for some $C>1$.  Namely, embeddedness
and the gradient estimate
(which applies because of stability) would
force  any nongraphical component
to spiral indefinitely, contradicting that $\Gamma$ is compact.  Thus it
is enough to find one component that is a graph.  This we be used below.

We will eventually show in Section \ref{s:areagrwth}
that (2) in Theorem
\ref{t:longstable}
does not happen; thus every component is a (single-valued) graph.
This will easily give Theorem \ref{t:remsing}.

\begin{figure}[htbp]
    \setlength{\captionindent}{20pt}
    \begin{minipage}[t]{0.5\textwidth}
    \centering\input{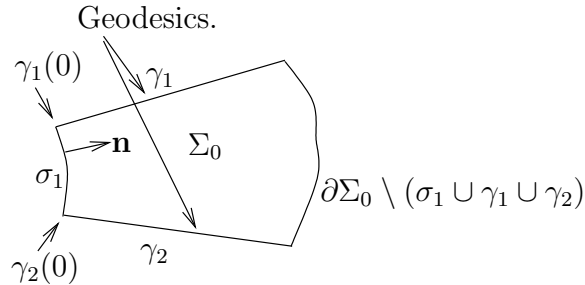}
    \caption{The subdomain $\Sigma_0 \subset \Sigma$ in Lemma
    \ref{l:stable} and below.}
    \label{f:sigma0}   \end{minipage}%
\end{figure}

\vskip2mm
See fig. \ref{f:sigma0}.
Throughout this section (except in Corollary
\ref{c:asmall}), $\Sigma \subset \RR^3$ will be an embedded
minimal planar domain (if the domain is stable, then we use
$\Gamma$ instead of $\Sigma$), $\Sigma_0\subset \Sigma$
a
subdomain and $\gamma_1$, $\gamma_2$, $\sigma_1 \subset \partial
\Sigma_0$ curves ($\gamma_1$, $\gamma_2$ geodesics) so
$\gamma_1\cup\gamma_2\cup\sigma_1$ is a simple curve and
$\gamma_i(0) \in \sigma_1$.  (By a geodesic we will mean
a curve with zero geodesic
curvature.  This definition of geodesic is needed when the curve intersect
the boundary of the surface.)  Below we will sometimes require
one or more of the following properties:\\
(A) $\dist_{\Sigma}(
\gamma_i(t) , \sigma_1 ) \geq t - C_0$ for $0 \leq t \leq
{\text{Length}} (\gamma_i)$.\\
(B) $\partial_{\nn} |x| \geq 0$ along
$\sigma_1$ (where $\nn$ is the inward normal to $
\partial \Sigma_0$).\\
(C) $\gamma_1\perp\sigma_1$, $\gamma_2\perp\sigma_1$ (i.e., angle $\pi/2$).\\
(D) $\dist_{\Sigma_0}( \sigma_1 , \partial \Sigma_0 \setminus
(\sigma_1 \cup \gamma_1 \cup \gamma_2) ) \geq \ell$ (thus $\ell \leq
{\text{Length}} (\gamma_i)$).\\

Note that if $\sigma_1 \subset \partial B_1$ (and $\Sigma_0$ is leaving $B_1$
along $\sigma_1$), then (B) is automatically
satisfied.

\vskip2mm The main component of the proof of Theorem
\ref{t:longstable} is Proposition \ref{c:bbl} below which shows
that certain stable sectors have
subsectors with small total curvature.  To show that we
will use an argument in the spirit of \cite{CM2}, \cite{CM4} to
get good curvature estimates for our nonstandard stable domains.
As in \cite{CM2}, \cite{CM4}, to estimate the total curvature we
show first an area bound.  That is (here $k_g$ is the geodesic
curvature of $\sigma_1$):

\begin{Lem} \label{l:stable}
Let $\Gamma_0 =\Gamma\subset \RR^3$ be stable and satisfy
(A) for $C_0=0$, (C), (D).
If $0 \leq \chi \leq 1$ is a function on
$\Gamma_0$ which vanishes on each $\gamma_i$, then for $1<R<\ell$
\begin{align}   \label{e:lstab}
    \Area (\an_{R} (\sigma_1 , \Gamma_0)) &\leq
C\,R^2 \,  \int_{\sigma_1}
   |k_g|
+ C \, R\,  {\text{Length}} (\sigma_1)\\
    &\quad
+C \, R^2 \,
    \left(  \int_{\an_{1}(\sigma_1 , \Gamma_0)} (1+ |A|^2 ) +
\int_{\an_{R}(\sigma_1 , \Gamma_0)} |\nabla \chi|^2 +
\int_{\an_{R}(\sigma_1 , \Gamma_0) \cap \{ \chi < 1 \} } |A|^2
\right) \, .\notag
\end{align}
\end{Lem}

\begin{proof}
Set $\an_{s,t} = \an_{s,t}(\sigma_1 , \Gamma_0)$ and $\dd =
\dist_{\Gamma} ( \sigma_1 , \cdot)$. Define a (radial) cut-off
function $\phi$ by
\begin{equation}
 \phi =
\begin{cases}
  \dd &\hbox{ on } \an_{1} \, ,    \\
  (R -\dd) / (R -1) &\hbox{ on } \an_{1,R}\, ,\\
  0 &\hbox{ otherwise }  \, . \\
\end{cases}
\end{equation}
By the stability inequality applied to $\phi\,\chi$ and using the
 inequality, $2ab\leq a^2+b^2$,
\begin{align}   \label{e:intcut}
    &\int_{\an_{1,R}}|A|^2\,[(R-\dd)/(R -1)]^2
    \leq \int |A|^2\,\phi^2 \leq 2 \int |\nabla \phi|^2
    +  2 \int_{\an_{R}}
 |\nabla \chi|^2 + \int_{ \an_{R}\cap \{ \chi < 1 \} }
    |A|^2   \notag \\
    &\quad \quad \leq  2\Area\,(\an_{1})+2(R -1)^{-2}\,
    \Area\,(\an_{1,R}) +  2 \int_{\an_{R}}
 |\nabla \chi|^2 + \int_{\an_{R} \cap \{ \chi < 1 \}}
    |A|^2   \, .
\end{align}
Set $K(s)= \int_{\an_{1,s}}|A|^2$.
By the coarea formula and integrating \eqr{e:intcut} by parts twice, we get
\begin{align} \label{e:intcut2}
    &2\,(R-1)^{-2}\int_1^{R} \int_1^t K(s)ds\,dt
    \leq 2/(R-1) \int_1^{R} K(s) (R-s)/(R -1)ds\notag \\
    &\quad \leq \int_{1}^{R}K'(s)  \,((R-s)/(R -1))^2 ds \\
    &\quad \quad \leq  2\Area\,(\an_{1})+2(R -1)^{-2}\,
    \Area\,(\an_{1,R}) +
    2 \int_{\an_{R}}
 |\nabla \chi|^2
     + \int_{\{ \chi < 1 \}}
    |A|^2    \notag\, .
\end{align}
Given $y \in \sigma_1$, let $\gamma_y:[0,r_y] \to \Gamma$ be the
(inward from $\partial \Gamma$) normal geodesic up to the
cut-locus of $\sigma_1$ (so $\dist_{\Gamma}(\sigma_1 , \gamma_y
(r_y) ) = r_y$) and $J_y$ the corresponding Jacobi field with $J_y
(0) = 1$ and $J_y'(0) = k_g (y)$.  Set $R_y = \min \{r_y , R\}$.
By the Jacobi equation,
\begin{equation}    \label{e:jy1}
    \int_0^{  R_y  } J_y(s) \, ds =
  R_y^2 \,    k_g (y) / 2
    + R_y      -
    \int_{0}^{R_y } \, \int_0^t  \, \int_0^{s}
    \K_{\Gamma}(\gamma_y (\tau)) \, J_y(\tau) \,
    d\tau \, ds \,dt \, .
\end{equation}
If $R_y<R$, then we extend $J_y(\tau)$,
$\K_y(\tau)=\K_{\Gamma}(\gamma_y(\tau))$
to functions $\tilde J_y$, $\tilde \K_y$ on $[0,R]$ by setting
$\tilde J_y=J_y$, $\tilde \K_y=\K_y$ on $[0,R_y]$ and
$\tilde J_y=\tilde \K_y=0$ otherwise.
If $R_y=R$, then we set $\tilde J_y=J_y$ and
$\tilde \K_y=\K_y$.
Since $\K_{\Gamma} = - |A|^2 /2$ (in particular, is $\leq 0$), by \eqr{e:jy1}
\begin{equation}    \label{e:jy10}
    \int_0^{  R_y  } J_y(s) \, ds \leq
  R^2 \,    |k_g (y)| / 2
    + R      -
    \int_{0}^{R } \, \int_0^t  \, \int_0^{s}
    \tilde \K_{y}(\tau) \, \tilde J_y(\tau) \,
    d\tau \, ds \,dt \, .
\end{equation}
Since
$K(s)=-2\int_{\sigma_1} \int_1^s
\tilde \K_y (\tau)\,\tilde J_y(\tau) \,d\tau\,dy$
(this uses (C)), integrating \eqr{e:jy10} over $\sigma_1$  gives
\begin{equation}   \label{e:intcut3}
    \Area\,(\an_{R})  \leq
\frac{R^2}{2} \,  \int_{\sigma_1}
   |k_g|
+  R\,  {\text{Length}} (\sigma_1)
+ \int_{1}^{R} \, \int_1^t  \frac{K(s)}{2} \, ds \,dt
    + \frac{R^2}{2} \int_{\an_{1}}|A|^2  \, .
\end{equation}
(Here we also used $\int_{0}^{R} \, \int_0^t  f(s) \, ds \, dt
\leq \int_{1}^{R} \, \int_1^t  [f(s) - f(1)]\, ds \, dt + R^2
\,f(1)$ for the nondecreasing function $f(t) = \int_{\an_t}|A|^2
\geq 0$.) Combining \eqr{e:intcut2} and \eqr{e:intcut3} gives
\eqr{e:lstab}.
\end{proof}

\begin{figure}[htbp]
    \setlength{\captionindent}{20pt}
    \begin{minipage}[t]{0.5\textwidth}
    \centering\input{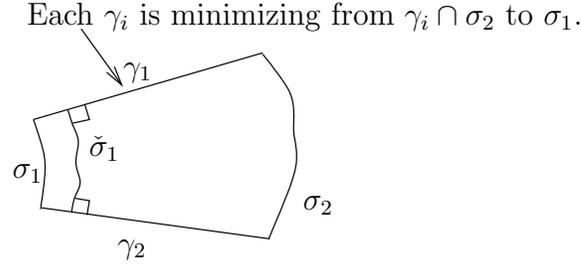}
    \caption{Lemma \ref{l:broken}: Connecting $\gamma_1$ and $\gamma_2$ by a curve
    $\check{\sigma}_1$ with  length and total curvature bounded.}
    \label{f:broken}   \end{minipage}%
\end{figure}

To apply Lemma \ref{l:stable}, we will need to replace a given
curve, in a minimal disk, by a curve lying within a fixed tubular
neighborhood of it and with length and total geodesic curvature
bounded in terms of the area of the tubular neighborhood.  This
is:

\begin{Lem}  \label{l:broken}
See fig. \ref{f:broken}.
If $\Sigma \subset \RR^3$ is an immersed minimal disk, $\partial
\Sigma = \gamma_1 \cup \gamma_2 \cup \sigma_1 \cup \sigma_2$, the
$\gamma_i$'s are geodesics with $2 \leq {\text{Length}}(\gamma_i)
= \dist_{\Sigma}(\sigma_2 \cap \gamma_i , \sigma_1)$, $1\leq
\dist_{\Sigma}(\sigma_1,\sigma_2)$, then there
exists a simple curve $ \check{\sigma}_1 \subset
\an_{1/64,1/4}(\sigma_1)$ connecting $\gamma_1$ to $\gamma_2$ and
with
\begin{equation}    \label{e:landtk}
 {\text{Length}}( \check{\sigma}_1 )
+ \int_{ \check{\sigma}_1 } |k_g| \leq C_1 \,
 (1+\Area\,(\an_{1/4}(\sigma_1))) \, .
\end{equation}
Moreover, we can choose $\check{\sigma}_1$ to intersect
$\gamma_i$ orthogonally and so
${\text{Length}}(\check{\gamma}_i)
= \dist_{\Sigma}(\sigma_2 \cap \gamma_i , \check{\sigma}_1)$, where
$\check{\gamma}_i$ denotes the component of
$\gamma_i\setminus \check{\sigma}_1$ which intersects $\sigma_2$.
\end{Lem}

\begin{proof}
We will do this in three steps.  First, we use
the coarea formula to find a level set of the
distance function with bounded length.
Local  replacement then gives a broken geodesic with
the same length bound and a bound on the number of
breaks.  Third, we find a simple subcurve and use
the Gauss-Bonnet theorem to control the number of breaks.

Set $\dd (\cdot) = \dist_{\Sigma} ( \sigma_1 , \cdot)$.  By the
coarea formula applied to (a regularization of) $\dd$, there
exists $d_0$ between $1/16$ and $3/32$ with ${\text{Length}} (\{
\dd = d_0 \}) \leq 32 \, \Area (\an_{1/8} (\sigma_1))$ and so that
$\{ \dd = d_0 \}$ is transverse. Since $\{ \dd = d_0 \}$ separates
$\sigma_1$ and $\sigma_2$, a component $\tilde \sigma$ of $\{ \dd
= d_0 \}$ goes from $\gamma_1$ to $\gamma_2$.

Parametrize $\tilde \sigma$ by arclength and let
$0=t_0<\cdots<t_n= \text{Length}(\tilde \sigma)$ be a subdivision
with $t_{i+1}-t_i \leq 1/32$ and $n \leq  32 \,
\text{Length}(\tilde \sigma) + 1$. Since $\cB_{1/32}(y)$ is a disk
for all $y\in \tilde \sigma$, it follows that we can replace
$\tilde \sigma$ with a broken geodesic $\tilde \sigma_1$ with
breaks at $\tilde \sigma(t_i)= \tilde \sigma_1 (t_i)$ and which is
homotopic to $\tilde \sigma$ in $\cT_{1/32}(\tilde \sigma)$.  We
can assume that $\tilde \sigma_1$ intersects the $\gamma_i$'s only
at its endpoints.

Let $[a,b]$ be a maximal interval so that $\tilde \sigma_1 |_{[a,b]}$
is simple.
We are done if $\tilde \sigma_1 |_{[a,b]} = \tilde \sigma_1$.
Otherwise, $\tilde \sigma_1 |_{[a,b]}$ bounds a disk in $\Sigma$ and
the Gauss-Bonnet
theorem implies that $\tilde \sigma_1 |_{(a,b)}$ contains a break.  Hence,
 replacing $\tilde \sigma_1$ by
$\tilde \sigma_1 \setminus \tilde \sigma_1 |_{(a,b)}$
gives a subcurve from $\gamma_1$ to $\gamma_2$ but does not
increase the number of breaks.  Repeating
this eventually gives a simple subcurve with the same bounds
for the length and the number of breaks.  Smoothing this at the breaks gives
the desired  $\check{\sigma}_1 $.

Finally, since $\gamma_i$ minimizes distance from $\gamma_i\cap\sigma_2$
to $\sigma_1$,  it follows easily by adding
segments in $\gamma_1 , \gamma_2$ to $\check{\sigma}_1$ and then
perturbing infinitesimally
near $\gamma_1$, $\gamma_2$ that
we can choose $\check{\sigma}_1$ to intersect $\gamma_i$ orthogonally
and so each
$\check{\gamma}_i$
minimizes distance back to $\check{\sigma}_1$; this gives at
most a bounded
contribution to the length
and total curvature.
\end{proof}

We will also need a version of Lemma \ref{l:broken}
where $\sigma$ is a
noncontractible curve (cf. lemma $1.21$ in \cite{CM4}):

\begin{Lem} \label{l:broken2}
Let $\Sigma \subset \RR^3$ be an immersed minimal planar domain and
$\sigma = B_1 \cap \partial \Sigma$  a simple closed curve with
$\dist_{\Sigma}(\sigma,\partial\Sigma \setminus \sigma) > 1$.
There
exists a simple noncontractible curve $ \check{\sigma} \subset
\an_{1/32,1/4}(\sigma)$ with
\begin{equation}    \label{e:landtk2}
 {\text{Length}}( \check{\sigma} )
+ \int_{ \check{\sigma} } |k_g| \leq C_1 \,
 (1+\Area\,(\an_{1/4}(\sigma))) \, .
\end{equation}
\end{Lem}

\begin{proof}
Following the first two steps of
the proof of Lemma \ref{l:broken} (with the
obvious modifications),  we get
a simple closed broken geodesic $\tilde \sigma_1$ which
is noncontractible with length and the number of breaks
$\leq C \, \Area\,(\an_{1/4}(\sigma))$.

As in the third step  of
the proof of Lemma \ref{l:broken}, let $\tilde \sigma_1 |_{[a,b]}$
be a maximal simple subcurve.  It follows that
$\tilde \sigma_1 |_{[a,b]}$ is closed (and has at most one
more break than $\tilde \sigma_1$).
If $\tilde \sigma_1 |_{[a,b]}$  is noncontractible, then we
are done.  Otherwise, if $\tilde \sigma_1 |_{[a,b]}$
bounds a disk, then we apply the Gauss-Bonnet theorem to see
that $\tilde \sigma_1 |_{(a,b)}$  contains a break
and proceed  as in  the proof of Lemma \ref{l:broken}.
\end{proof}

In Proposition \ref{c:bbl} below, we will also need a
lower bound for the area growth of tubular neighborhoods of a
curve.  To get such a bound, we need that the curve is not
completely ``crumpled up.'' This will follow by requiring that
$(t+C_0) \,(t+1)\leq  \delta \, \Area (\an_{1}(\sigma_1))$.

\begin{Lem}  \label{l:areagwth}
Let $\Sigma_0 = \Sigma$ satisfy  (A), (B) and (D). If $\sigma_1
\subset B_1$, $1 \leq s < t \leq \ell$ and $(t+C_0) \,(t+1)\leq
\delta \, \Area (\an_{1}(\sigma_1))$, then
\begin{equation}    \label{e:stda}
 (t+1)^{2\delta-2}  \, \Area (\an_{t}(\sigma_1 ))    \geq
(s+1)^{2\delta-2}  \,  \Area\,(\an_{s}(\sigma_1))  \, .
\end{equation}
\end{Lem}

\begin{proof}
Set $\an_t = \an_t (\sigma_1)$ and $L(s) = \int_{
\partial \an_s \setminus \partial \Sigma} 1$. By minimality,
Stokes' theorem, (A), (B) and
$\dist_{\Sigma}(\sigma_1,x)+1\geq |x|$, we get that
\begin{equation}
4\,\Area\, (\an_{s})
=\int_{\an_{s}}\Delta \,|x|^2
\leq 2\,(s+1)\,L(s)+4\, (s+C_0) \,(s+1)\, .
\end{equation}
By the coarea formula,
$(\Area \, (\an_{s}))' = L(s)$ for almost every $s$.
Hence,
 for almost every $s$ with
$\dist_{\Sigma} (\sigma_1 , \sigma_2) \geq s\geq 1$
\begin{equation}
(\log\,\Area \,(\an_{s}))'\geq \frac{2}{s+1} -\frac{2\,(s+C_0)}{\Area
(\an_{s})}\geq \frac{2\,(1-\delta)}{s+1} \, .
\end{equation}
Since $\Area (\an_s)$ is monotone, a standard argument then gives
\eqr{e:stda}.
\end{proof}

\begin{Rem} \label{r:spc}
In the special case of Lemma \ref{l:areagwth} where $\Sigma$ is an
annulus with $\partial \Sigma=\sigma_1 \cup \sigma_2$, i.e., where
$\gamma_i=\emptyset$ and $\sigma_1$, $\sigma_2$ are closed,
 the proof
simplifies in an obvious way and $\delta$ can be chosen to be
zero.
\end{Rem}

We are now ready to apply Lemma \ref{l:stable} and use the
logarithmic cut-off trick to show that certain stable sectors have
small curvature:

\begin{Pro} \label{c:bbl}
Let $\Gamma_0 \subset \Gamma \subset \RR^3$ satisfy
(A) (with $C_0 = 0$), (B),
(D), and $\dist_{\Gamma} (\Gamma_0 , \partial \Gamma) > 1/4$.
 Suppose that $\Gamma$ is
stable,
$\omega > 2$,  $\ell > R_0 > \omega^2$, and $\sigma_1 \subset B_1$.
If $\Gamma_0$ is a disk and
$4\,R_0^2 \,(R_0+1)\leq \Area
(\an_{1}(\sigma_1 , \Gamma_0))$, then for $\omega^2 \leq t \leq R_0$
\begin{align}   \label{e:quada}
\Area\,(\an_{2} (\sigma_1,\Gamma_0))\,t^{2}/C &\leq  \Area (
\an_{\omega , t}(\sigma_1,\Gamma_0) )
    \leq C \, \Area\,(\an_{2} (\sigma_1,\Gamma_0))\,t^2
   \, ,    \\
    \int_{\an_{\omega , R_0 / \omega}(\sigma_1 , \Gamma_0)} |A|^2
&\leq C\, R_0+\frac{C}{\log \omega}  \Area\,(\an_{2} (\sigma_1,\Gamma_0))
\, .   \label{e:totcurv}
\end{align}
\end{Pro}

\begin{proof}
Define a  function $\chi$ on  $\Gamma_0$ by
\begin{equation}
 \chi =
\begin{cases}
  2 \, \dist_{\Gamma}(\gamma_1 \cup \gamma_2 ,
    \cdot )
  &\hbox{ on } \cT_{1/2}(\gamma_1 \cup \gamma_2) \, ,    \\
  1 &\hbox{ otherwise }  \, .\\
\end{cases}
\end{equation}
We will use $\chi$ to cut-off on the sides $\gamma_1 ,
\gamma_2$.
 Using \cite{Sc},
\cite{CM2},  and $\dist_{\Gamma} (\Gamma_0 , \partial \Gamma) > 1/4$,
\begin{align}
    \int_{\an_{2} (\sigma_1, \Gamma_0) } (1+|A|^2 )
&\leq C_1 \, \Area (
\an_{2}(\sigma_1,\Gamma_0)  ) \, , \label{e:thischi}\\
 2 \int_{\an_{R_0}(\sigma_1 , \Gamma_0)}
 |\nabla \chi|^2
+ \int_{ \an_{R_0}(\sigma_1 , \Gamma_0) \cap \{ \chi < 1 \} }
|A|^2 &\leq C_1 \, R_0 \leq C_1  \, \Area (
\an_{1}(\sigma_1,\Gamma_0) ) \, .\label{e:thischi0}
\end{align}
Since $\sigma_1 \subset \partial \Gamma_0$
satisfies (A) with $C_0=0$ and (D),
Lemma \ref{l:broken} gives
a simple curve  $\check{\sigma}_1$  (and $\check{\gamma}_1$,
$\check{\gamma}_2$)
satisfying  (A) with $C_0=0$, (C),  (D), and \eqr{e:landtk};
let $\check{\Gamma}_0 \subset \Gamma_0$ be the component
of $\Gamma_0\setminus \check{\sigma}_1$ containing $\sigma_2$.
By the triangle inequality,
\begin{equation}  \label{e:btti}
    \an_{t} (\check{\sigma}_1,\Gamma_0) \subset \an_{t+1/4} (\sigma_1,\Gamma_0)
\subset \an_{t+1/4} (\check{\sigma}_1, \check{\Gamma}_0)  \cup
(\Gamma_0 \setminus \check{\Gamma}_0) \, .
\end{equation}
Note that $\Gamma_0 \setminus \check{\Gamma}_0$ is a disk with
boundary $\sigma_1 \cup \check{\sigma}_1 \cup (\gamma_1 \setminus
\check{\gamma}_1) \cup (\gamma_2 \setminus \check{\gamma}_2)$.
Hence, by minimality, Stokes' theorem,  (B),  $|x|\leq 5/4$ on
$\partial (\Gamma_0 \setminus \check{\Gamma}_0)$,
and \eqr{e:landtk},
\begin{equation}    \label{e:bdleftover}
    4 \, \Area (\Gamma_0 \setminus \check{\Gamma}_0)
    = \int_{\Gamma_0 \setminus \check{\Gamma}_0} \Delta |x|^2
\leq 2 \, \int_{\check{\sigma}_1 \cup (\gamma_1 \setminus
\check{\gamma}_1) \cup (\gamma_2 \setminus \check{\gamma}_2)}
    |x|  \leq
    C_1' \, \Area (
\an_{1}(\sigma_1,\Gamma_0) )  \, .
\end{equation}
 Inserting \eqr{e:thischi}, \eqr{e:thischi0} into Lemma
\ref{l:stable} applied to $\check{\sigma}_1$ and using
\eqr{e:landtk}, \eqr{e:btti},  \eqr{e:bdleftover}
 gives for $2 \leq t
\leq R_0$
\begin{equation}  \label{e:areab}
    \Area ( \an_{t}(\sigma_1,\Gamma_0) )
    \leq C_2 \, \Area\,(\an_{2} (\sigma_1,\Gamma_0)
) \, t^2 \, ,
\end{equation}
giving the second inequality in \eqr{e:quada}.
 Set $\an_{t} = \an_{t}
(\sigma_1,\Gamma_0)$ (define $\an_{s,t}$ similarly) and set
$L(t) = \int_{\partial
\an_t \setminus \partial \Gamma_0} 1$.
By \eqr{e:areab}, the coarea formula, and integration
by parts,
\begin{align}    \label{e:log1}
\int_{R_0 / \omega}^{R_0} L(t)\,t^{-2}dt&=
\left[ \Area (\an_{R_0/\omega, t})\,t^{-2} \right]_{R_0 / \omega}^{R_0}
+2\, \int_{R_0 / \omega}^{R_0}
\Area (\an_{R_0/\omega, t})\,t^{-3} \, dt \notag\\
&\leq  C_2\,(1+2\log\,\omega )\,\Area\,(\an_{2})\leq C_3\,\log \omega\,
\Area\,(\an_{2})\, ,\\
\int_{1}^{\omega} L(t)\,t^{-2}dt &\leq
\Area (\an_{1,\omega}) \, \omega^{-2}+2\, \int_1^{\omega}
\Area (\an_{1,t})\,t^{-3} \, dt  \leq
C_3\,\log\, \omega
\,\Area\,(\an_{2}) \, .\label{e:log2}
\end{align}
Define a (radial) cut-off function $\eta$ by
\begin{equation}
 \eta =
\begin{cases}
  \log \dist_{\Gamma_0} ( \sigma_1  , \cdot) /
\log \omega &\hbox{ on } \an_{1,\omega} \, ,    \\
  1 &\hbox{ on } \an_{\omega,R_0 / \omega}\, ,\\
  \left[ \log R_0 - \log \dist_{\Gamma_0} ( \sigma_1  , \cdot)
    \right] /
    \log \omega &\hbox{ on } \an_{R_0 / \omega,R_0} \, .\\
\end{cases}
\end{equation}
Using the bounds \eqr{e:log1}
and \eqr{e:log2}, we get
\begin{align}   \label{e:tpopt}
    \int |\nabla \eta|^2 &= \int_{ \an_{1,\omega} } |\nabla \eta|^2
    + \int_{ \an_{R_0 / \omega, R_0} } |\nabla \eta|^2   \\
    &\leq  \frac{1}{ (\log \omega)^2 } \int_{1}^{\omega}
\frac{L(t)}{t^2}dt
    + \frac{1}{ (\log  \omega)^2 }
\int_{R_0 / \omega}^{R_0} \frac{L(t)}{t^2}dt
        \leq \frac{C_3 \, \Area\,(\an_{2})}{\log \omega} \, . \notag
\end{align}
   Substituting $\eta \, \chi$ into the stability inequality,
we get using \eqr{e:thischi0} and \eqr{e:tpopt}
\begin{equation}    \label{e:stabsays}
    \int_{ \an_{\omega,R_0 / \omega} } |A|^2
\leq \int_{ \an_{R_0} \cap \{ \chi < 1 \} } |A|^2 +2
\int_{\an_{R_0}}
 |\nabla \chi|^2
+ 2 \,  \int |\nabla \eta|^2 \leq
    C_1\,R_0+\frac{2\,C_3 \,\Area\,(\an_{2})}{\log \omega}     \, .
\end{equation}
Finally, Lemma
\ref{l:areagwth} (and \eqr{e:areab} for $t=\omega$)
gives the first inequality in \eqr{e:quada}.
\end{proof}

We will prove Theorem \ref{t:longstable} by considering two
separate cases depending on the area of $\an_1(\sigma)$.  When
$\Area ( \an_1(\sigma))$ is small, the next corollary
will show that (1) of
Theorem \ref{t:longstable} holds.  On the
other hand, when $\Area ( \an_1(\sigma))$ is large, we will show in
the next section, using Corollary
\ref{c:multigr} below, that (2) of Theorem \ref{t:longstable} holds.

\begin{Cor}     \label{c:asmall}
Given $C_a$, there exists $\Omega_a > 4$ so: Let $\Gamma \subset
\RR^3$ be a stable embedded minimal planar domain,
$\sigma = B_1 \cap \partial \Gamma$ connected, and
$\dist_{\Gamma}(\sigma,\partial\Gamma \setminus \sigma)> R$.
If  $R > \Omega_a^2$ and $\Area (\an_1(\sigma)) \leq C_a$, then
$\Gamma$ contains a graph $\Gamma_g$ (after a rotation) over
$D_{R/\Omega_a} \setminus
D_{\Omega_a}$ with gradient $\leq 1$ and $\dist_{\Gamma} (\sigma , \Gamma_g)
\leq 2 \, \Omega_a$.
\end{Cor}

\begin{proof}
Lemma \ref{l:broken2} gives a simple closed noncontractible curve
$ \check{\sigma} \subset
\an_{1/32,1/4}(\sigma)$ and with
 ${\text{Length}}( \check{\sigma} ) + \int_{ \check{\sigma} } |k_g| \leq C_1 \,
 [\Area\,(\an_{1}(\sigma))+1] $.  Since $\Gamma$ is a planar domain,
$\check{\sigma}$ separates in $\Gamma$; let $\check{\Gamma}$ be
the component of  $\Gamma \setminus \check{\sigma}$  which does not
contain $\sigma$.
By Lemma \ref{l:stable} (which applies with $\chi \equiv 1$ since
$\gamma_1=\gamma_2=\emptyset$), we get for $1 \leq t \leq R$
\begin{equation}    \label{e:qag}
    \Area (  \an_t(\check{\sigma}, \check{\Gamma})  )
\leq C \, (C_a+1) \, t^2 \, .
\end{equation}
Given $\Omega > 4$, by
\eqr{e:qag}
 and the logarithmic cut-off trick
in the stability inequality (cf. \eqr{e:stabsays}), we get
that $\int_{\an_{\Omega/2 ,2 R/\Omega} (  \check{\sigma}, \check{\Gamma}  )}
|A|^2 \leq C_2 \, (C_a+1)  / \log \Omega$.   Combining this with \eqr{e:qag}
and the Cauchy-Schwarz inequality gives for $\Omega / 2 \leq t \leq R/\Omega$
\begin{equation}    \label{e:qag2}
    \int_{\an_{t,2t} (\check{\sigma}, \check{\Gamma})}
|A| \leq \left( \Area (\an_{2t}( \check{\sigma}, \check{\Gamma}))
\int_{\an_{\Omega/2 ,2 R/\Omega} (
    \check{\sigma}, \check{\Gamma})}
|A|^2 \right)^{1/2}  \leq \frac{C_3 \, (C_a+1)  \, t}{(\log \Omega)^{1/2}} \, .
\end{equation}
Applying the coarea formula on $\an_{t,2t}$ for $t=\Omega/2 ,
R/\Omega$, \eqr{e:qag2} gives a (possibly  disconnected)
planar domain $\Gamma_0 \subset \an_{\Omega/2 ,2 R/\Omega}
(\check{\sigma}, \check{\Gamma})$
with  $\an_{\Omega , R/\Omega} (\check{\sigma}, \check{\Gamma})
\subset \Gamma_0$,
$\partial \Gamma_0 = \cup_{i=1}^n \sigma_i$, and
\begin{equation}    \label{e:qah}
    \sum_{i=1}^n \int_{\sigma_i} |A|
\leq \frac{C_3 \, (C_a+1) }{(\log \Omega)^{1/2}} \, .
\end{equation}
We now fix $\Omega = \Omega (C_a) > 4$ so that
$C_2 \, (C_a+1)  / \log \Omega < \pi$ and
$C_3 \, (C_a+1)  \, (\log \Omega)^{-1/2} < 1/4$.  Using $\int_{\Gamma_0}
|A|^2 < \pi$,
\eqr{e:qah} (which implies equation $(1.13)$ in \cite{CM7} with
$\epsilon = 1/4$), and
that the Gauss map is conformal,  proposition $1.12$ of \cite{CM7}
 implies that  on each component $\Gamma_0^k$ of $\Gamma_0$ we get
 $\nn_{\Gamma}(\Gamma_0^k)\subset \cB_{1/2}(a_k)$, where
$a_k \in \SS^2$ ($\nn_{\Gamma}$ is the unit normal to $\Gamma$).
Hence,
$\Gamma_0$ is a (possibly multi-valued) graph.  Since $\Gamma$ is
embedded,  the corollary now follows easily (cf.
 lemma $1.10$ in \cite{CM4}).
\end{proof}

We construct next from curves
$\sigma_1$, $\gamma_1$, $\gamma_2$ in a stable surface the
desired multi-valued graph.  (The existence of the
curves
$\sigma_1$, $\gamma_1$, $\gamma_2$ will be established
in the next section.)  First we need the following two
lemmas:

\begin{Lem}     \label{l:suc}
Given $C_1 , \epsilon_0 > 0$, there exists $\epsilon_1 > 0$ so
if $\cB_1 \subset \Sigma$ is minimal, $\sup_{\cB_{1/2}} |A|^2 \leq \epsilon_1$,
and $\sup_{\cB_{1}} |A|^2 \leq C_1$, then
$\sup_{\cB_{3/4}} |A|^2 \leq \epsilon_0$.
\end{Lem}

\begin{proof}
Suppose not; it follows that there is a sequence
$\Sigma_j$ of minimal surfaces with
$\sup_{\cB_{1/2}} |A|^2 \leq 1/j$,
 $\sup_{\cB_{1}} |A|^2 \leq C_1$, and
$\sup_{\cB_{3/4}} |A|^2 >  \epsilon_0 > 0$.  The uniform bound
$\sup_{\cB_{1}} |A|^2 \leq C_1$ (and standard elliptic estimates)
gives a subsequence which converges in $C^{2,\alpha}$
to a limit $\Sigma_{\infty}$.
It follows that $\Sigma_{\infty}$ is minimal, $|A|^2 = 0$ on $\cB_{1/2}$,
and
$\sup_{\cB_{3/4}} |A|^2 \geq  \epsilon_0 > 0$.  By unique continuation,
$\Sigma_{\infty}$ is flat contradicting that
$\sup_{\cB_{3/4}} |A|^2 \geq  \epsilon_0 > 0$.
\end{proof}

The next lemma will be applied both when $\Gamma$ is an annulus
and when $\Gamma$ has boundary on the sides.  When
$\Gamma$ is an annulus,
the condition \eqr{e:tilded} will be trivially satisfied and
it will be possible for $\Gamma$ to contain a graph instead of
a multi-valued graph.

\begin{figure}[htbp]
    \setlength{\captionindent}{20pt}
    \begin{minipage}[t]{0.5\textwidth}
    \centering\input{pl10.pstex_t}
    \caption{The proof of
    Lemma \ref{l:techlem}: Repeatedly applying Lemma \ref{l:suc} along
    chains of balls builds out a ``flat'' region in $\Gamma$.}
    \label{f:techlem}   \end{minipage}%
\end{figure}

\begin{Lem}         \label{l:techlem}
Given $N , S_0 >4$ , $\epsilon > 0$, there exist $C_b > 1$,
 $\delta > 0$  so: Let
$\Gamma\subset \RR^3$ be a stable embedded minimal surface and
$\sigma = B_1 \cap \partial \Gamma$. If $\gamma : [0,S_0] \to
\Gamma$ is a geodesic with $\dist_{\Gamma} (\gamma (t) , \sigma) =
t$ for $0\leq t \leq S_0$,
  $\sup_{\cB_{S_0/16}(\gamma (S_0))} |A|^2
\leq \delta \, S_0^{-2}$, and for $0\leq t \leq S_0$
\begin{equation}   \label{e:tilded}
    \dist_{\Gamma \setminus \an_{t/8} (\sigma) } ( \gamma (t) ,
    \partial \Gamma ) \geq
        C_b \, t \, ,
\end{equation}
then (after a rotation of $\RR^3$) $\Gamma$ contains either an
$N$-valued graph $\Gamma_{2, S_0/2}^{-N\pi,N\pi}$, or a graph
$\Gamma_{2, S_0/2}$,
 with gradient $\leq \epsilon$,
$|A| \leq \epsilon/r$,  and
$\gamma (4) \in  \Gamma_{2, 5}^{-\pi,\pi}$ (or in
$\Gamma_{2, 5}$).
\end{Lem}

\begin{proof}
Combining estimates for stable surfaces of \cite{Sc}, \cite{CM2}
and \eqr{e:tilded}, gives for $0\leq t\leq S_0$
\begin{equation}   \label{e:pcb0}
    \sup_{\cB_{t/2} (\gamma (t))} |A|
    \leq  C_0 \, t^{-1}  \, .
\end{equation}
Fix $\delta_0> 0$  to be chosen small depending on $S_0$. Using
\eqr{e:pcb0} and repeatedly applying Lemma \ref{l:suc} along a
chain of balls with centers in $\gamma$, see fig. \ref{f:techlem},
there exists $\delta_1 =
\delta_1 (S_0 , \delta_0 , C_0) > 0$ so that if $\delta \leq
\delta_1$, then
 for $1\leq t\leq S_0$
\begin{equation}   \label{e:pcb1}
    \sup_{\cB_{t/32} (\gamma (t))} |A|
    \leq   \delta_0 \, t^{-1}  \, .
\end{equation}
 Since
$\gamma$ is a geodesic in $\Gamma$,  \eqr{e:pcb1} gives the bound
$k_g^{\RR^3} (t) \leq \delta_0 \, t^{-1}$ for the geodesic
curvature of $\gamma$ in $\RR^3$. It follows that for $1 \leq t
\leq S_0$
\begin{equation}    \label{e:intgga}
   |\nn_{\Gamma} (\gamma(t)) - \nn_{\Gamma} (\gamma(1))|
+ |\gamma'(t) - \gamma'(1)|
    \leq  2 \delta_0 \, \int_1^{S_0} \frac{ds}{s}
    \leq 2 \delta_0 \, \log S_0 \, ;
\end{equation}
i.e.,  $\gamma$ is $C^1$-close to a straight line segment
in $\RR^3$ and
$\nn_{\Gamma}$ is almost constant on $\gamma$.
Rotate so that $\gamma'(1) = (1,0,0)$ (i.e.,
so $\gamma'(1)$  points in the $x_1$-direction).
 For $\delta_0 > 0$ small, \eqr{e:intgga} (and $\gamma (0) \in B_1$)
implies that for $1\leq t \leq S_0$
\begin{equation}    \label{e:proper}
    3t/4 - 2 \leq x_1 (\gamma(t)) \leq 1 + t  \, .
\end{equation}

We will now argue as in \eqr{e:pcb0}--\eqr{e:pcb1} to extend the
region where $\Gamma$ is graphical, this time using balls centered
on cylinders (i.e., building out the multi-valued graph in the
$\theta$ direction). Suppose now that $4\leq s \leq S_0/2$ and
$y_{0,s} = \{ x_1^2 + x_2^2 = s^2 \} \cap \gamma$. Using
\eqr{e:pcb1}, $\cB_{C_2 s}(y_{0,s})$ is a graph with gradient
$\leq C_2' \, \delta_0$ over $\nn_{\Gamma} (y_{0,s})$. In
particular, also using \eqr{e:intgga}, $\partial \cB_{C_2
s}(y_{0,s})$  contains a point $y_{1,s} \in \{ x_1^2 + x_2^2 = s^2
\}$.  Using Lemma \ref{l:suc}, we can therefore repeat this to
find $y_{2,s}$, etc.  It follows from \eqr{e:tilded} that we can
continue this until $\Gamma$ either closes up (giving a graph) or
we have the desired
 $N$-valued graph $\Gamma_{2,S_0/2}^{-N\pi,N\pi}$
with gradient $\leq \epsilon$, $|A| \leq
\epsilon/r$,  and which contains
$\gamma (4)$.
\end{proof}

In the next corollary $\Gamma\subset B_{2R}\subset \RR^3$ will be
a stable embedded minimal annulus with $\partial \Gamma\subset
B_{1/4}\cup \partial B_{2R}$ where $B_1 \cap \partial \Gamma$ is
connected and suppose $\Gamma_0 \subset \Gamma$ is a disk
satisfying (A) for $C_0=0$, (B), (D). Let $\sigma=B_1\cap \partial
\Gamma$ so $\sigma_1\subset \sigma$ and $\sigma$ is a simple
closed curve. Assume also that the following
strengthening of (A) holds:\\
(A') $\dist_{\Gamma}(\gamma_i(t),\sigma)=t$ for
$0\leq t\leq \text{Length}(\gamma_i)$.

\begin{figure}[htbp]
    \setlength{\captionindent}{20pt}
       \begin{minipage}[t]{0.5\textwidth}
    \centering\input{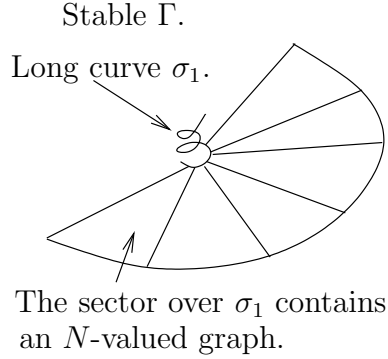}
    \caption{Corollary \ref{c:multigr}:  A stable tubular neighborhood of
    a long curve $\sigma_1$ contains an $N$-valued
    graph $\Gamma_{\omega_0, R/\omega_0}^{-N\pi,N\pi}$.}
    \label{f:multigr}   \end{minipage}
\end{figure}

\begin{Cor}  \label{c:multigr}
See fig. \ref{f:multigr}.
Given $N , \epsilon > 0$, there exist $\omega_0, R_0 > 1$ so if
$\Gamma$, $\Gamma_0$ are as above, and
$\Area (\cT_{1}(\sigma_{1}))\geq 4\,R_0^2 \,(R_0+1)$,
then (after a rotation of $\RR^3$) $\Gamma$
contains an $N$-valued graph $\Gamma_{\omega_0, R/\omega_0}^{-N\pi,N\pi}$
with gradient $\leq \epsilon$,
$|A| \leq \epsilon/r$,  and
\begin{equation}
\dist_{\Gamma} ( z_1 , \Gamma_{\omega_0, \omega_0}^{0,0})
< 2 \, \omega_0+ C_1\,\Area (\an_{1}(\sigma_1 , \Gamma))\, .
\end{equation}
\end{Cor}

\begin{proof}
 Proposition \ref{c:bbl} gives $C$ so that for $\omega \leq t \leq R_0 / \omega$
 (where $\omega > 4$ and $R_0^2 > 2 \omega^2$)
\begin{align}   \label{e:quadaz}
\Area\,(\an_{2} (\sigma_1,\Gamma_0))\, t^{2}/C &\leq  \Area (
\an_{\omega , t}(\sigma_1,\Gamma_0) )
       \, ,    \\
    \int_{\an_{\omega , R_0 / \omega}(\sigma_1 , \Gamma_0)} |A|^2
&\leq \frac{C}{\log \omega}  \Area\,(\an_{2} (\sigma_1,\Gamma_0))
\, .   \label{e:totcurvz}
\end{align}
(Here we also used $\Area (\cT_{1}(\sigma_{1}))\geq 4\,R_0^2
\,(R_0+1)$ in \eqr{e:totcurvz}.) Set $S=\omega$. Choose a maximal
disjoint collection of balls $\cB_{S/4}(y_1) , \dots ,
\cB_{S/4}(y_n)$ with centers in $\an_{S,2S}(\sigma_1,\Gamma_0)$.
Since $\Gamma$ is annulus without boundary on the sides and $R_0 >
5S/2$, it follows from (A') that $\cB_{S/2}(y_j) \cap
\partial \Gamma = \emptyset$; we use this twice.  First,
 since $\an_{S,2S}(\sigma_1,\Gamma_0)$
is contained in the union of the double balls and $\pi (S/4)^2
\leq \Area (\cB_{S/4}(y_j)) \leq \Area (\cB_{S/2}(y_j)) \leq C \pi
\, S^2$ by stability (see \cite{CM2}), we  have $n \geq C \,
S^{-2} \, \Area\,(\an_{S,2S}(\sigma_1,\Gamma_0))$. Second, again
by stability, \cite{CM2}, $\int_{\an_{S/4}(\gamma_1 \cup \gamma_2)
\cap \an_{S/2,3S}(\sigma_1) } |A|^2 \leq C$.
 Combining this
with \eqr{e:quadaz} and
\eqr{e:totcurvz},  we can find $j$ so that
$\int_{\cB_{S/4}(y_j)} |A|^2 < C / \log \omega$,
and therefore, by the mean value inequality,
$\sup_{\cB_{S/8}(y_j)} |A|^2 < C \, S^{-2} / \log \omega$.

Let $\gamma:[0,\ell] \to \Gamma$ be a minimal geodesic from
$y_j$ to $\sigma_{1}$;
note that $S \leq \ell \leq 2S$.
Using that the sides $\gamma_{1} , \gamma_{2}$ are minimizing
(i.e., (A')), it follows that
$\gamma \subset \Gamma_{0}$.  Furthermore, since $\Gamma$ is an annulus, (A')
implies that
\begin{equation}   \label{e:tildeda}
    \dist_{\Gamma \setminus \an_{1} (\sigma) } ( \gamma (\ell) ,
    \partial \Gamma ) \geq  R / 2
    \, ,
\end{equation}
In particular,  given $\omega_1 , N_1 > 1$ and and $\epsilon_1 > 0$,
there exists $\omega$ (and hence $R_0$) large  so we
can apply Lemma \ref{l:techlem} to
get either a graph $\Gamma_{S/\omega_1 , S/2}$
or an initial multi-valued graph
$\Gamma_{S / \omega_1 , S/2}^{-N_1\pi,N_1\pi}$ with
gradient $\leq \epsilon_1$,
$|A| \leq \epsilon_1 /r$,  and
$\gamma (4S/\omega_1) \in  \Gamma_{2S/\omega_1, 5S/\omega_1}^{-\pi,\pi}$.
However, since $\Area (\cT_{1}(\sigma_{1}))\geq 4\,R_0^2 \,(R_0+1)$,
 $\Gamma$ cannot contain  a graph $\Gamma_{S/\omega_1 , S/2}$.

Using theorem II.0.21 of \cite{CM3}, we will next extend
$\Gamma_{S / \omega_1, S/2}^{-N_1\pi,N_1\pi}$ to the desired
$N$-valued graph $\Gamma_{\omega_0, R/\omega_0}^{-N\pi,N\pi}$.
Namely, let $P$ be the vertical plane $\{ x_1 = 2S /\omega_1 \}$.
We claim first that each component of $P\cap \Gamma$ goes off to
$\partial B_{2R}$. To see this, note that by the maximum
principle, any closed curve in $P \cap \Gamma$ would be homologous
to the interior boundary of $\Gamma$ and together these two curves
would span an annulus in $\Gamma$ violating the convex hull
property (using the multi-valued graph in $\Gamma$ to connect this
annulus to  $\{ x_1 = -S /\omega_1 \}$). It follows that two of
these nodal curves connect the multi-valued graph out to $\partial
B_{2R}$, giving a curve $\eta$ in $\Gamma$ with both endpoints in
$\partial B_{2R}$. One component of $\Gamma \setminus \eta$ is a
stable disk which is forced to spiral initially. Therefore, by
theorem II.0.21 of \cite{CM3}, this extends to the desired
multi-valued graph.
\end{proof}

\section{The minimizing geodesics and the proof of Theorem \ref{t:longstable}}

In Proposition \ref{p:constructgamma} and Corollary
\ref{c:constructgamma} below,
we will construct the minimizing geodesics $\gamma_1 , \gamma_2$
needed for Corollary \ref{c:multigr}.
To do that we will first need the following lemmas and corollaries
(here $\cT_t$ is the closed tubular neighborhood and $\cT_t^{\circ}$
is the open):

\begin{figure}[htbp]
    \setlength{\captionindent}{20pt}
    \begin{minipage}[t]{0.5\textwidth}
    \centering\input{pl14.pstex_t}
    \caption{The set $E$ in Lemma \ref{l:regl1}.}
    \label{f:reg0}   \end{minipage}\begin{minipage}[t]{0.5\textwidth}
    \centering\input{pl15.pstex_t}
    \caption{In an annulus $\Sigma$ with $\partial \Sigma = \sigma_1 \cup \sigma_2$,
    given geodesics $\gamma_1$, $\gamma_2$ and a curve $\sigma_3 \subset \sigma_1$
    connecting $\gamma_1(0)$ and $\gamma_2(0)$, Lemma \ref{l:regl1}
    finds a disk $\Sigma_4$ with $\partial \Sigma_4 = \sigma_3 \cup \sigma_4 \cup
    \gamma_1 \cup \gamma_2$ where each point in $\sigma_4$ is almost distance $\ell$
    from $\sigma_1$.}
    \label{f:reg}   \end{minipage}
\end{figure}

\begin{Lem}   \label{l:regl1}
See fig. \ref{f:reg0} and
\ref{f:reg}.
   Let $\Sigma$ be an annulus with
$\partial \Sigma=\sigma_1\cup\sigma_2$, where
$\dist_{\Sigma}(\sigma_1,\sigma_2)>\ell+\epsilon$ for $\ell$, $\epsilon>0$
and let $E$
be the connected component of
$\Sigma\setminus \cT_{\ell}(\sigma_1)$ containing $\sigma_2$.
Let $\gamma_1$,
$\gamma_2$ be geodesics with
\begin{equation}   \label{e:rege1}
\text{$\gamma_i:[0,\ell]\to \Sigma$,
$\dist_{\Sigma}(\gamma_i(t),\sigma_1)=t$ for $0\leq t\leq \ell$,
and $\gamma_i(\ell)\in \overline{E}$}\, .
\end{equation}
If $\sigma_3\subset \sigma_1$ is a segment
connecting $\gamma_1(0)$ and $\gamma_2(0)$, then
there exists a curve
$\sigma_4\subset \cT_{\epsilon}^{\circ} (E)
\cap \cT_{\epsilon}^{\circ}(\Sigma\setminus E)$ connecting $\gamma_1(\ell)$
and $\gamma_2(\ell)$ and so
$\sigma_3\cup\sigma_4\cup \gamma_1\cup\gamma_2$ bounds a disk
$\Sigma_4$.  Moreover,
$\sigma_4\subset \cT_{\ell+\epsilon}(\sigma_1)\setminus
\cT_{\ell-\epsilon}(\sigma_1)$.
\end{Lem}

\begin{proof}
First, note that $\gamma_1(\ell)$,
$\gamma_2(\ell)\in \overline{E}\cap\overline{\Sigma\setminus E}$ and by
definition $E$, hence $\cT_{\epsilon}^{\circ}(E)$, is connected.  Moreover,
if $x\in \Sigma\setminus E$ and $\gamma:[0,\ell_{\gamma}]\to \Sigma$
is a geodesic with
$\gamma (\ell_{\gamma})=x$ and $\dist_{\Sigma}(\gamma (t),\sigma_1)=t$ for
$0\leq t \leq \ell_{\gamma}$, then $\gamma\cap E=\emptyset$.
Hence, also $\Sigma\setminus E$
and $\cT_{\epsilon}^{\circ}(\Sigma\setminus E)$ are connected.  Since
$\sigma_1 \subset \cT_{\epsilon}^{\circ}(\Sigma\setminus E)$ and
$\sigma_2 \subset \cT_{\epsilon}^{\circ}( E)$,
applying van Kampen's theorem to
$\Sigma = \cT_{\epsilon}^{\circ}(\Sigma\setminus E) \cup
\cT_{\epsilon}^{\circ}( E)$ gives that
$\cT_{\epsilon}^{\circ}(E)
\cap \cT_{\epsilon}^{\circ}(\Sigma\setminus E)$ is path connected and has
fundamental group $\ZZ$ which injects into $\pi_1(\Sigma)$.
In particular, we get simple curves
$\sigma_{4,1} , \sigma_{4,2} \subset \cT_{\epsilon}^{\circ}(E)
\cap \cT_{\epsilon}^{\circ}(\Sigma\setminus E)$
connecting $\gamma_1(\ell)$
to $\gamma_2(\ell)$ and so
$\sigma_{4,1} \cup \sigma_{4,2}$ is homologous to $\sigma_1$.
Fix $\Sigma_0 \subset \Sigma$ with $\partial \Sigma_0 =
\sigma_1 \cup (\sigma_{4,1} \cup \sigma_{4,2})$.
The curve
$\sigma_3 \cup \gamma_1\cup\gamma_2$ divides $\Sigma_0$ into
two components, one of which is a disk with $\sigma_3$, $\gamma_1$,
$\gamma_2$, and either $\sigma_{4,1}$ or $\sigma_{4,2}$ in its boundary.

Finally, since $\sigma_4\subset \cT_{\epsilon}^{\circ}(E)$ it follows that
$\sigma_4\subset \Sigma\setminus \cT_{\ell-\epsilon}(\sigma_1)$.
Likewise it follows from that
$\sigma_4\subset \cT_{\epsilon}^{\circ}(E)
\cap\cT_{\epsilon}^{\circ}(\Sigma\setminus E)$ and the triangle inequality
that $\sigma_4\subset \cT_{\ell+\epsilon}(\sigma_1)$.
\end{proof}

\begin{figure}[htbp]
    \setlength{\captionindent}{20pt}
    \begin{minipage}[t]{0.5\textwidth}
    \centering\input{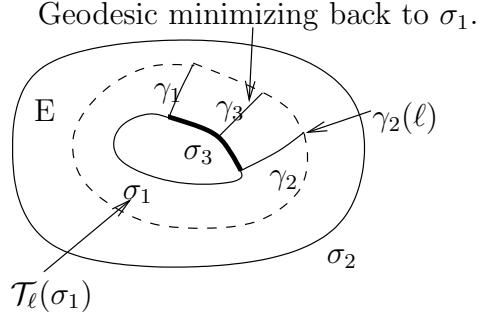}
    \caption{Corollary \ref{c:regc1}: Finding a  geodesic
    $\gamma_3$ satisfying \eqr{e:rege1} between
    two other geodesics $\gamma_1$, $\gamma_2$.}
    \label{f:regc}   \end{minipage}
\end{figure}

\begin{Cor}   \label{c:regc1}
See fig. \ref{f:regc}.
Let $\Sigma$, $E$, $\sigma_1$, $\sigma_2$, $\sigma_3$, $\gamma_1$,
$\gamma_2$
 be as in Lemma \ref{l:regl1}.   If $\gamma_1(\ell)\ne \gamma_2(\ell)$, then
there exists a geodesic
$\gamma_3$ different from $\gamma_1,\,\gamma_2$, intersecting
$\sigma_3$, and satisfying \eqr{e:rege1}.
\end{Cor}

\begin{proof}
Let $\eta\subset \Sigma\setminus (\gamma_1\cup\gamma_2)$
be a simple curve from $\sigma_3$ to $\sigma_2$ and so
$\eta\cap \sigma_1\subset \eta\cap \sigma_3$ is one point.  Fix $\mu>0$ with
$3\mu<\dist_{\Sigma}(\gamma_1\cup\gamma_2,\eta)$. For
$\epsilon>0$ small (in particular,
$\epsilon<\dist_{\Sigma}(\sigma_1,\sigma_2)-\ell$),
let $\sigma_{\epsilon,4}$, $\Sigma_{\epsilon,4}$
be given by Lemma
\ref{l:regl1}.
Let $\eta_{\epsilon}$
be the component of $\eta\cap \Sigma_{\epsilon,4}$
intersecting $\sigma_3$ and let
$\gamma_{\epsilon,3}:[0,\ell_{\epsilon}]\to \Sigma$ be a geodesic
with $\gamma_{\epsilon,3}(\ell_{\epsilon})
=\partial \eta_{\epsilon}\setminus \sigma_1$ and
$\dist_{\Sigma}(\gamma_{\epsilon,3}(t),\sigma_1)=t$ for
$0\leq t\leq \ell_{\epsilon}$. Since
$\sigma_{\epsilon,4}\subset \cT_{\ell+\epsilon}(\sigma_1)
\setminus \cT_{\ell-\epsilon}(\sigma_1)$, then
$\ell-\epsilon< \ell_{\epsilon}\leq \ell+\epsilon$.  Moreover, by the
triangle inequality, if $\epsilon<\mu$, then
$\cB_{\mu}(\gamma_k(\ell))\cap \gamma_{\epsilon,3}=\emptyset$
for $k=1,2$, hence
$\cB_{\mu}(\gamma_k(\ell))\cap (\eta\cup \gamma_{\epsilon,3})=\emptyset$
for $k=1,2$ and $\epsilon<\mu$.  We claim that
\begin{equation}  \label{e:keyeq}
\text{$\eta_{\epsilon}\cup \gamma_{\epsilon,3}\subset
\cT_{\delta}(\Sigma\setminus E)$ where $\delta\to 0$ as $\epsilon\to 0$}\, .
\end{equation}
Suppose that \eqr{e:keyeq} fails; it follows that there exists a sequence
$\epsilon_i\to 0$ and $x_i\in\eta_{\epsilon_i}$ with $x_i\to x$, where
$x_i,x\in\Sigma\setminus \cT^{\circ}_{\delta}(\Sigma\setminus E)\subset E$
for some $\delta>0$.
Since $E$ is open and connected, there exists a curve $\nu\subset E$ from $x$
to $\sigma_2$ and so
$\nu\subset E\setminus \cT_{\delta_0}(\Sigma\setminus E)$ for
some $\delta_0=\delta_0(x)>0$.  For $i$ sufficiently large, we
can extend $\nu$ to a curve
$\nu_i\subset E\setminus \cT_{\delta_0/2}(\Sigma\setminus E)$
from $x_i$ to $\sigma_2$.  However,
$\gamma_1\cup\gamma_2\cup\sigma_{\epsilon_i,4}
\subset \cT_{\epsilon_i}(\Sigma\setminus E)$
separates $x_i$ from $\sigma_2$ which is a contradiction for $i$
sufficiently large.  Hence, \eqr{e:keyeq} holds.

Pick a sequence $\epsilon_i>0$ with $\epsilon_i\to 0$. After
passing to a subsequence, we can assume that
 $\gamma_{\epsilon_j,3}\to \gamma_3$.  It is clear that
$\gamma_{3}:[0,\ell]\to \Sigma$ is a geodesic with
$\gamma_{3}(0)\in \sigma_{1} \setminus \{ \gamma_1(0), \gamma_2
(0) \}$, $\dist_{\Sigma}(\gamma_{3}(t),\sigma_1)=t$ for
$0\leq t\leq \ell$, and $\gamma_3(\ell)\in \overline{E}$. It remains
to show that $\gamma_{3}(0)\in \sigma_{3}$.

If $\gamma_{3}(0) \notin \sigma_{3}$, then $\dist_{\Sigma}
(\gamma_3(0) , \sigma_3 ) > 0$ (since $\gamma_{3}(0)\in \sigma_{1}
\setminus \{ \gamma_1(0), \gamma_2 (0) \}$) and therefore
$\dist_{\Sigma} (\gamma_{\epsilon_j,3}(0) , \sigma_3 ) > 0$ for
$j$  large. It follows that $\eta_{\epsilon_j}\cup \gamma_{\epsilon_j,3}$
divides
$\Sigma$ into two components
$\Sigma_{\epsilon_j,1}$,
$\Sigma_{\epsilon_j,2}$ with
$\gamma_1 (\ell) \in \Sigma_{\epsilon_j,1}$ and
$\gamma_2 (\ell) \in \Sigma_{\epsilon_j,2}$. (That $\gamma_1(\ell)$,
$\gamma_2(\ell)$ are in different components follows from that
$\gamma_{\epsilon,3}
\cap \gamma_1=\emptyset=\gamma_{\epsilon,3}\cap \gamma_2$ by
the triangle inequality.) After
possibly switching $\gamma_1$ and $\gamma_2$ (and going to a subsequence),
we can assume that
$\sigma_2 \subset \Sigma_{\epsilon_j,2}$.   Note that,
$\cB_{\mu}(\gamma_1(\ell))\subset \Sigma_{\epsilon_j,1}$ since we showed
above that
$\cB_{\mu}(\gamma_1(\ell))\cap (\eta\cup\gamma_{\epsilon,3})=\emptyset$.
We will use this
to contradict that $\gamma_1(\ell)\in \overline {E}$.  Namely, let
$x\in \cB_{\mu/2}(\gamma_1(\ell))\cap E$ (note  that such an $x$ exists
since $\gamma_1(\ell)\in \overline{E}$).  Since $E$ is open and connected,
then there exists a curve
$\nu\subset E\setminus \cT_{\delta_0}(\Sigma\setminus E)$ for some
sufficiently small $\delta_0=\delta_0(x)>0$ which connect $x$ and
$\sigma_2$.  This contradicts \eqr{e:keyeq} for $j$ sufficiently large
since $\eta_{\epsilon_j}\cup\gamma_{\epsilon_j,3}$ separate
$\Sigma_{\epsilon_j,1}$ and $\sigma_2$.
\end{proof}

\begin{Lem} \label{l:isop}
If $\Sigma \subset \RR^3$ is an immersed minimal surface with
$\partial \Sigma = \gamma_1 \cup \gamma_2 \cup \sigma$ where
$\sigma \subset B_1$, $\partial_{\nn}|x|\geq 0$ on $\sigma$ ($\nn$
is the inward normal to $\partial \Sigma$), and
$\gamma_1,\,\gamma_2$ have length $\leq \ell$, then
\begin{equation}    \label{e:abda}
    \Area (\an_1(\sigma)) \leq  4\,\ell \,(\ell+1) \, .
\end{equation}
\end{Lem}

\begin{proof}
By minimality, Stokes' theorem, $\partial_{\nn}|x|\geq 0$ on
$\sigma$, and $|x|\leq \ell +1$ on $\gamma_i$,
\begin{equation}
    4 \, \Area (\Sigma) = \int_{\Sigma} \Delta |x|^2
\leq 2 \, \int_{\gamma_1 \cup \gamma_2}
    |x| \, |\partial_{\nn} |x|| \leq
    4 \, (\ell+1) \, \ell \, .
\end{equation}
\end{proof}

In what follows, if $\sigma\subset \partial \Sigma$
is a simple curve, $\nn$ is the inward normal to $\sigma$,
$\tilde \sigma$ is a segment of $\sigma$, then
(see fig. \ref{f:inwhat})
\begin{equation}        \label{e:inwhat}
\cT_s(\tilde\sigma,\nn)
=\{\exp_{\tilde\sigma (t)}(\tau\,\nn(t))\,
|\,\dist_{\Sigma}(\exp_{\tilde\sigma (t)}(\tau\,\nn(t)),\sigma)
=\tau\leq s\}\, .
\end{equation}

\begin{figure}[htbp]
    \setlength{\captionindent}{20pt}
    \begin{minipage}[t]{0.5\textwidth}
    \centering\input{pl17.pstex_t}
    \caption{The region $\cT_s(\tilde\sigma,\nn)$ in \eqr{e:inwhat}.}
    \label{f:inwhat}   \end{minipage}%
    \begin{minipage}[t]{0.5\textwidth}
    \centering\input{pl18.pstex_t}
    \caption{Proposition \ref{p:constructgamma}:
    Finding a geodesic $\gamma \subset \Sigma$ which minimizes
    back to the curve $\sigma$.}
    \label{f:congamma}   \end{minipage}
\end{figure}

\begin{Pro} \label{p:constructgamma}
See fig. \ref{f:congamma}.
Let $\Sigma \subset \RR^3$ be an immersed minimal annulus,
$\sigma\subset B_1\cap \partial \Sigma$ a simple closed curve
with $\dist_{\Sigma}(\sigma,\partial\Sigma\setminus \sigma)> \ell$,
$\partial_{\nn}|x|\geq 0$ on $\sigma$, and let $E$
be as in Lemma \ref{l:regl1}.
If $\tilde \sigma$ is a segment of $\sigma$ and
$\Area (\cT_1(\tilde \sigma,\nn))> 4\,\ell\, (\ell+1)$,
then there exists a geodesic
$\gamma:[0,\ell ]\to \Sigma$ with $\dist_{\Sigma}(\gamma(t),\sigma)=t$
for $0\leq t\leq \ell$ and $\gamma (0)\in \tilde \sigma$,
$\gamma (\ell)\in \overline{E}$.
\end{Pro}

\begin{proof}
Suppose that there is no such geodesic $\gamma$.
Let $B$ be the set of geodesics satisfying \eqr{e:rege1} for $\sigma_1=\sigma$.
It follows easily that $A=\{\gamma_0(0)\,|\,\text{ $\gamma_0\in B$}\}$
is a closed subset of $\sigma\setminus \tilde \sigma$
containing more than two points.  Let $\hat \sigma$ be the
connected component of $\sigma\setminus A$ containing
$\tilde \sigma$ (note that $\hat \sigma$ is open) and
let $\partial \hat \sigma=\{\gamma_1(0),\gamma_2(0)\}$
where $\gamma_1$, $\gamma_2$ are the corresponding minimizing
geodesics of lengths $\ell$.

By Corollary \ref{c:regc1}, $\gamma_1(\ell)=\gamma_2(\ell)$.
In fact,
there exists a subset $\hat \Sigma$ of $\Sigma$ with
$\partial \hat \Sigma=\gamma_1\cup \gamma_2\cup \hat \sigma$.  Since
$\Area (\cT_1(\hat \sigma,\tilde \Sigma))\geq
\Area (\cT_1(\tilde \sigma,\tilde \Sigma))
=\Area (\cT_1(\tilde \sigma,\nn))> 4 \,\ell\, (\ell+1)$, it follows from
Lemma \ref{l:isop} that $A\cap \tilde \sigma\ne \emptyset$
which is the desired contradiction and the proposition follows.
\end{proof}

\begin{Cor} \label{c:constructgamma}
Let $\Sigma \subset \RR^3$ be an immersed minimal annulus,
 $\sigma\subset B_1\cap \partial \Sigma$
is a simple closed curve with
$\dist_{\Sigma}(\sigma,\partial\Sigma\setminus \sigma)> \ell \geq
1$, $\partial_{\nn}|x|\geq 0$ on $\sigma$, and $\Area
(\cT_1(\sigma,\Sigma))> 12\,\ell^2 \, (\ell+1)$. For each
$z_1\in\sigma$ there is a segment $z_1\in\sigma_1\subset \sigma$
and geodesics $\gamma_1$, $\gamma_2:[0,\ell ]\to \Sigma$ with
$\{\gamma_1(0),\gamma_2(0)\}=\partial \sigma_1$,
\begin{equation}  \label{e:consteq1}
\text{$\dist_{\Sigma}(\gamma_i(t),\sigma)=t$ for $0\leq t\leq
\ell$}\, .
\end{equation}
Moreover, for all $\epsilon>0$ a disk $\Sigma_0\subset\Sigma$ has $\sigma_1$,
$\gamma_i\subset
\partial \Sigma_0$, $\dist_{\Sigma_0} (\partial \Sigma_0\setminus
\sigma_1\cup \gamma_1\cup\gamma_2)> \ell-\epsilon$,
\begin{equation}  \label{e:consteq2}
15\,\ell^2 \, (\ell+1)>\Area (\cT_1(\sigma_1,\Sigma_0))> 4\,\ell^2 \,
(\ell+1)\, .
\end{equation}
\end{Cor}

\begin{proof}
Let $\sigma_{z_1}^1$, $\sigma_{z_1}^2$, $\sigma_{z_1}^3$ be three
consecutive (disjoint) subsegments of $\sigma$ with
$z_1\in\sigma_{z_1}^2$ being the ``middle one'' and so for each $i$
\begin{equation}   \label{e:consteq3}
5\,\ell^2 \, (\ell+1) > \Area
(\cT_1(\sigma_{z_1}^i,\nn))>4\,\ell^2 \, (\ell+1)\, .
\end{equation}
By Proposition \ref{p:constructgamma} applied to both
$\sigma_{z_1}^1$ and $\sigma_{z_1}^3$, we get geodesics
$\gamma_1$, $\gamma_2:[0,\ell]\to \Sigma$ satisfying
\eqr{e:consteq1} and with $\gamma_1(0)\in \sigma_{z_1}^1$,
$\gamma_2(0)\in \sigma_{z_1}^3$, $\gamma_i(\ell)\in \overline{E}$
(where $E$ is the connected component of $\Sigma\setminus
\cT_{\ell}(\sigma)$ containing $\sigma_2$).
 Let $\sigma_1$ be the segment of $\sigma$ between $\gamma_1(0)$
and $\gamma_2(0)$ and containing $\sigma_{z_1}^2$. By Lemma
\ref{l:regl1} there is a disk $\Sigma_0\subset\Sigma$ with
$\sigma_1$, $\gamma_1$, $\gamma_2\subset
\partial \Sigma_0$, and $\dist_{\Sigma_0} (\partial \Sigma_0\setminus
\sigma_1\cup \gamma_1\cup\gamma_2)> \ell-\epsilon$.  We need to show
\eqr{e:consteq2}.  Since $\sigma_{z_1}^2\subset \sigma_1$ the
lower bound in \eqr{e:consteq2} follows easily from
\eqr{e:consteq3}. To see the upper bound, observe that if $x\in
\cT_1(\sigma_1,\Sigma_0)$, then clearly
$\dist_{\Sigma}(x,\sigma)=\dist_{\Sigma_0}(x,\sigma_1)$ and hence
\begin{equation}   \label{e:consteq4}
\cT_1(\sigma_1,\Sigma_0)\subset
\cup_{i=1,2,3}\cT_1(\sigma_{z_1}^i,\nn)\, .
\end{equation}
 From \eqr{e:consteq3} and \eqr{e:consteq4}, the upper bound in
\eqr{e:consteq2} follows.
\end{proof}

\begin{proof}
(of Theorem \ref{t:longstable}).
Given $N$, $\epsilon$,
let $\omega_0$, $R_0$ be given by Corollary \ref{c:multigr}.  Set
$\sigma = \partial B_1 \cap \partial \Gamma_1 (\partial)$ and note
that $\partial_{\nn} |x| \geq 0$ and $B_1 \cap \partial \Gamma$ is
connected since $\Gamma$ is an annulus.  Suppose that
 $\dist_{\Sigma}(\sigma,\partial \Gamma\setminus \sigma)>R_0$.
By Corollary \ref{c:asmall}, (1)
holds if  $\Area(\an_1 (\sigma))\leq 12\,R_0^2 \,(R_0+1)$.
(Recall that if one component of
$B_{R/\omega} \cap \Gamma \setminus B_{\omega}$ contains a graph over
$D_{R/(2\omega)} \setminus
D_{2 \omega}$ with gradient $\leq 1$, then every component of
$B_{R/(C\omega)} \cap \Gamma \setminus
B_{C\omega}$ is a graph for some $C>1$.)

On the other hand if $\Area(\an_1 (\sigma))> 12\,R_0^2 \,(R_0+1)$,
then it follows from Corollary \ref{c:multigr}
together with Corollary \ref{c:constructgamma} that (2) holds.
\end{proof}

Using that the curvature of a $2$-valued embedded minimal graph
decays faster than quadratically (this was shown in \cite{CM8}),
we show next (which will be needed in the next section)
that such $2$-valued graphs contain minimal geodesics
close to the
radial curve  $\theta=0$. (In particular, there is such
a geodesic which do not spiral.)  In this corollary,
$\Gamma_{\lambda \omega}  (\partial)$  denotes the component of
$B_{\lambda \omega} \cap \Gamma$ containing $B_{\lambda \omega} \cap
\partial \Gamma$.

\begin{Cor}  \label{c:longstable}
There exists $\lambda > 1$ so: If $\Gamma$ is as in Theorem
\ref{t:longstable}, $\Gamma_{\omega,R/\omega}^{-3\pi,3\pi}$ is as
in (2) of that theorem (with $N \geq 3$, $\tau \leq 1$), and $R >
\lambda \omega^2$, then there exists a geodesic
$\gamma:[0,\ell]\to \Gamma_{\omega , R^{1/2}}^{-\pi,\pi}$ with
$\gamma (0)\in \partial B_{\lambda \omega}$,
$\dist_{\Gamma}(\gamma (t), \Gamma_{\lambda \omega }(\partial)
)=t$ for $0 \leq t \leq \ell$, and $\ell \geq R^{1/2}/4$.
\end{Cor}

\begin{proof}
Fix $\lambda > 1$ large to be chosen. Set $\dd = \dist_{
\Gamma_{\omega  , R^{1/2}  }^{-2\pi,2\pi}} ( \Gamma_{\omega
,\omega  }^{-2\pi,2\pi}  , \cdot)$ and let $\Gamma_{\omega
 ,R/\omega}^{-3\pi,3\pi}$
be the graph of $u$. By Corollary 1.14 of \cite{CM8}, on
$S_{\omega , R^{1/2}}^{-2\pi,2\pi}$ we have $\rho \, |\Hess_u|
\leq C' \, (\rho/\omega)^{-5/12}$. Hence,  on $\Gamma_{\omega
 ,  R^{1/2}  }^{-2\pi,2\pi}$
\begin{equation}    \label{e:ceedd}
   \dd \, |A| \leq C \, \omega^{5/12} \, \dd^{-5/12} \, .
\end{equation}
 Let $\gamma:[0,\ell]\to \Gamma$ be a minimizing geodesic in $\Gamma$
 from (the point) $\Gamma_{ R^{1/2}/3  , R^{1/2}/3  }^{0,0}$
to $\Gamma_{\lambda \omega}(\partial)$, so $\dist_{\Gamma}(\gamma
(t), \Gamma_{\lambda \omega} (\partial) )=t$ for $0\leq t\leq
\ell$. In particular, $\gamma(0) \in
\partial \Gamma_{\lambda \omega}(\partial)$ and $\gamma(\ell)
= \Gamma_{ R^{1/2}/3  ,  R^{1/2} /3 }^{0,0}$. Using the radial
curve $\Gamma_{\omega  ,  R^{1/2} /3 }^{0,0}$ as a comparison (and
that $\tau \leq 1$), $ {\text{Length}}(\gamma) = \ell \leq
R^{1/2}/2$. Let $\tilde{\gamma}$ be the maximal segment of
$\gamma$ in $\Gamma_{\omega  ,  R^{1/2} }^{-\pi,\pi}$ containing
$\gamma(\ell)$.  Since $\tilde{\gamma}$ is a geodesic in $\Gamma$,
\eqr{e:ceedd} gives the bound $k_g^{\RR^3} (t) \leq C\,
\omega^{5/12} \, t^{-1-5/12}$ for the geodesic curvature of
$\tilde{\gamma}$ in $\RR^3$. It follows that for $\lambda \omega
\leq t \leq \ell$
\begin{equation}    \label{e:intgg}
    |\tilde{\gamma}'(t) - \gamma'(\ell)| \leq C
    \, \omega^{5/12} \,  \int_{\lambda \omega}^{\infty}
    s^{-17/12} \, ds \leq  12\, C \, \lambda^{-5/12} / 5 \, ;
\end{equation}
i.e., $\tilde{\gamma}$ is $C^1$-close to a straight line segment
in $\RR^3$.   For $\lambda$ large, \eqr{e:intgg} implies that
either $\tilde{\gamma} \subset \Gamma_{\omega , R^{1/2}/2
}^{-3\pi/4,3\pi/4}$ or $\gamma$ leaves $B_{ R^{1/2}}$; the latter
is impossible since ${\text{Length}}(\gamma) \leq  R^{1/2} /2 $.
We conclude that $\tilde{\gamma} \subset \Gamma_{\omega
 , R^{1/2}/2 }^{-3\pi/4,3\pi/4}$. In particular,
$\tilde{\gamma} = \gamma$ and the corollary follows.
\end{proof}

\section{Area growth of stable sectors
and the proof of Theorem \ref{t:remsing}}
 \label{s:areagrwth}

In this section, we show that case (2) in Theorem
\ref{t:longstable} does not happen and thus Theorem
\ref{t:remsing} follows easily. To do that, we first prove upper
and lower bounds for the area of a stable sector over a curve
$\sigma_1$ if the  sides $\gamma_1 ,\gamma_2$ of the sector are
contained in multi-valued graphs $\Sigma_1 , \Sigma_2$. By
\cite{CM8}, the number of sheets of each $\Sigma_i$ grows at least
like $\log^2 \rho$, giving the lower area bound $\rho^2 \, \log^2
\rho$ when the $\Sigma_i$'s are disjoint. We use this growing
number of sheets to construct a function $\chi$ with small energy
and which vanishes on the sides $\gamma_1 ,\gamma_2$. Inserting
$\chi$ in Lemma \ref{l:stable} gives the upper area bound $\rho^2
(C + \log \log \rho)$ (where $C= C(\sigma_1)$). If $\rho$ is large
depending on $C$, then these bounds are contradictory and hence
the $\Sigma_i$'s cannot be disjoint.

We will use several times that given $\alpha > 0$, Proposition
II.2.12 of \cite{CM3} gives $N_g > 0$ so if $u$ satisfies the
minimal surface equation on $S_{\e^{-N_g}, \e^{N_g} \, R}^{- N_g,
2 \pi + N_g}$ with $|\nabla u| \leq 1$, and $w<0$ (where $w$ is
the separation), then $\rho \, |\Hess_u| + \rho \, |\nabla  w
|/|w| \leq \alpha$ on $S_{1, R}^{0 , 2\pi }$. Theorem 3.36 of
\cite{CM7} then yields $|\nabla u - \nabla u(1,0)| \leq C \alpha$.
We can therefore assume (after rotating so $\nabla u(1,0) = 0$)
that
\begin{equation} \label{e:wantit}
    |\nabla u| + \rho \, |\Hess_u| +  4 \, \rho \, |\nabla  w |/|w| +
    \rho^2 \, |\Hess_w  |/|w|
      \leq \epsilon < 1/(2\pi) \, .
\end{equation}
The bound on $|\Hess_w  |$ follows from the other bounds and
standard elliptic theory.

The next lemma shows that an embedded
 multi-valued minimal graph in a concave cone (intersected with  cylindrical
 shells; see fig. \ref{f:cone})
\begin{equation}
    \cC_{\Lambda ,  R}(h)
    =\{ x \, | \, (x_3 - h)^2 \leq  \Lambda^2 \, (x_1^2 + x_2^2) \,
    ,\, 1/4 \leq x_1^2 + x_2^2 \leq R^2 \}
\end{equation}
 has at least $\log^2 \rho$ many
 sheets.   Note that the axis of the cone $\cC_{\Lambda ,  R}(h)$
 is the $x_3$-axis and the vertex
 is $(0,0,h)$.
 We will only need $\log \rho$ sheets for most of what
 follows, except for the lower bound for the area given
 in Corollary  \ref{c:area} below.

\begin{figure}[htbp]
    \setlength{\captionindent}{20pt}
    \begin{minipage}[t]{0.5\textwidth}
    \centering\input{pl19.pstex_t}
    \caption{The truncated cone  $\cC_{\Lambda ,  R}(h)$.}
    \label{f:cone}   \end{minipage}%
    \begin{minipage}[t]{0.5\textwidth}
    \centering\input{pl20.pstex_t}
    \caption{Lemma \ref{l:cone1}: It takes at least $\approx \log^2 \rho$
    rotations for a multi-valued graph to spiral out of the cone
    $\cC_{8\, \pi \, \epsilon, R}(0)$.}
    \label{f:lcone1}   \end{minipage}
\end{figure}

\begin{Lem} \label{l:cone1}
See fig. \ref{f:lcone1}.
Given $\epsilon > 0$, there exist
$0< C_1 < 1$ and $C_2$ so: Let  $\Sigma\subset \cC_{8
\, \pi \, \epsilon, R}(0)$ with $\partial \Sigma \subset
\partial \cC_{8 \, \pi \, \epsilon,R}(0)$ be a
minimal multi-valued graph of $u$
with $w<0$ and $u(1,0)=0$.  If the domain of $u$ contains
$S_{1/2,R}^{-2\pi,2\pi}$ and $u$ satisfies \eqr{e:wantit}, then
$\Sigma$ contains a (multi-valued) graph over
\begin{equation}   \label{e:cgd}
\{ (\rho , \theta) \, | \, |\theta| \leq C_1 \, \log^2 \rho + \pi
, \, 1 \leq \rho \leq R^{3/4} / 2 \}
 \end{equation} with
\begin{equation}    \label{e:cee}
    \rho^2 \, |A|^2 \leq C_2 \, \rho^{-5/18} \, .
\end{equation}
\end{Lem}

\begin{proof}
Corollary 1.14 of \cite{CM8} gives
on $S_{1,R^{3/4}}^{-\pi,\pi}$ that
\begin{equation}    \label{e:cee1}
    \rho^2 \, |\Hess_u (\rho , \theta) |^2 \leq C \, \rho^{-5/18} \, ,
\end{equation}
directly giving \eqr{e:cee} for $|\theta| \leq  \pi$.
By corollary 5.7 of \cite{CM8}, $\Sigma$ contains a (multi-valued)
graph over $\{ (\rho , \theta) \, | \, c^2 \, |\theta| \leq \log^2
\rho , \, 1 \leq \rho \leq R^{3/4} \}$  with $|u(\rho ,2\pi n) -
u(\rho,0)| \leq \rho^{\epsilon}$ for $n \in \ZZ$,  $2 \, \pi \,
c^2 \, |n| \leq \log^2 \rho$.
Applying the Harnack inequality and elliptic estimates to the
function $w_n(\rho,\theta) = u(\rho ,2\pi n + \theta) -
u(\rho,\theta)$ (cf. (1.17) of \cite{CM8}), we get
\begin{equation}    \label{e:bess}
    \rho \,
|\nabla u(\rho ,2\pi n) - \nabla      u(\rho,0)|
    + \rho^2 \,
 |\Hess_u (\rho ,2\pi n) - \Hess_u(\rho,0)|
\leq C' \, \rho^{\epsilon}  \, .
\end{equation}
Combining  \eqr{e:cee1} and \eqr{e:bess} then easily gives
\eqr{e:cee} in general.
\end{proof}

  We first define  a function $0 \leq \chi \leq 1$  on $\cP$
(the universal cover
  of $\CC \setminus \{ 0 \}$) which
is $0$ on $S_{3/4,\infty}^{-\pi , \pi}$,   $1$ on $\{ \rho <
R^{3/4}/2\} \setminus ( S_{1/2,R}^{-2\pi , 2\pi} \cup$
\eqr{e:cgd}), and so $|\nabla_{\cP} \chi|^2$ is of the order
$(\rho\, \log \rho)^{-2}$ for $\rho$ large. Namely, set
\begin{equation}    \label{e:chionp}
  \chi (\rho , \theta) =
\begin{cases}
3- 4\rho
&{\text{ for }} 1/2 \leq \rho < 3/4 , \,  |\theta| \leq  \pi  \, ,\\
1-(C_1 - |\theta| + \pi)(4\rho - 2)/C_1
&{\text{ for }} 1/2 \leq \rho < 3/4 , \, \pi \leq  |\theta|
\leq  C_1 + \pi  \, ,\\
0&{\text{ for }}  |\theta| \leq \pi , \,  3/4 \leq \rho  \, ,\\
(|\theta| - \pi) /C_1
&{\text{ for }} 3/4 \leq \rho < \e , \,  \pi \leq |\theta|
\leq C_1+ \pi  \, ,\\
  (|\theta| - \pi) /(C_1 \, \log \rho )&{\text{ for }} \e \leq \rho , \, \pi
\leq | \theta| \leq C_1 \,
\log \rho + \pi \, ,\\ 1 &\hbox{ otherwise} \, .
\\
\end{cases}
\end{equation}
Using \eqr{e:chionp}, define $\chi$ on a (multi-valued) graph over
a domain containing $S_{1/2,R}^{-2\pi,2\pi} \cup$ \eqr{e:cgd}
 in the obvious way.
 Note that if $\Sigma$ is as in Lemma \ref{l:cone1}, then
$1-\chi$ is one on the central
sheet $\Sigma_{3/4,R}^{-\pi , \pi}$
and vanishes before $\Sigma$ leaves the cone on the top, bottom,
or inside.

\begin{Cor} \label{c:bdchi}
Given $\epsilon > 0$,
there exists $C_3$ so if $\Sigma , u $ are as in Lemma
\ref{l:cone1},
then $\chi =0$ over
$S_{3/4,R}^{-\pi , \pi}$, $\chi = 1$ on
$\{ x_1^2 + x_2^2 \leq R^{3/2}/4 \} \cap \partial
\Sigma$, and
for $\e < t \leq R^{3/4}/2$
\begin{equation}    \label{e:nocurv}
    \int_{\{ \chi < 1 , \, x_1^2 + x_2^2 \leq t^2 \} } |A|^2
    + \int_{\{ x_1^2 + x_2^2 \leq t^2 \} \cap \Sigma} |\nabla \chi|^2
     \leq C_3 \, ( 1 + \log \log t )  \,  .
\end{equation}
\end{Cor}

\begin{proof}
Clearly, $\chi =0$ over
$S_{3/4,R}^{-\pi , \pi}$. By Lemma \ref{l:cone1},
 $\chi = 1$ on $\{ x_1^2 + x_2^2 \leq R^{3/2}/4 \} \cap \partial
\Sigma$.
To get \eqr{e:nocurv},
first consider $\chi$ as a function downstairs on $\cP$. On $\{
\rho \leq \e \}$, $|\nabla_{\cP} \chi| \leq C_0$ and $\{
|\nabla_{\cP} \chi| \ne 0 \} \subset \{ |\theta| \leq C_1 + \pi ,
\, 1/2 \leq \rho\}$. Similarly, on $\{ \e \leq \rho \}$,
$|\partial_{\theta} \, \chi|/\rho \leq 1/ (C_1 \, \rho \, \log
\rho)$ and $|\partial_{\rho} \, \chi| \leq 1 / (\rho \, \log
\rho)$, so that $ |\nabla_{\cP} \chi|^2   \leq 2 \, (C_1 \,  \rho
\, \log \rho)^{-2}$ and $\{ |\nabla_{\cP} \chi| \ne 0 \} \subset
\{ \pi \leq |\theta | \leq C_1 \log \rho + \pi  \}$.
Therefore,
since $\Sigma$ is a graph with gradient $\leq 1$, it follows
easily that
\begin{equation}    \label{e:echib}
   \int_{\{ x_1^2 + x_2^2 \leq t^2 \} \cap \Sigma} |\nabla \chi|^2
\leq
C_0' + \frac{12}{ C_1} \, \int_{\e}^t
    \frac{ds}{s\, \log s}  = C_0' + \frac{12 \log \log t}{C_1}
\, .
\end{equation}
Similarly, using \eqr{e:cee} gives
\begin{equation}    \label{e:areabb}
   \int_{\{ \chi < 1 , \, x_1^2 + x_2^2 \leq t^2 \} } |A|^2
\leq C+
 4  \, C_2
\int_{\e}^{\infty} (\pi + C_1 \, \log s)
     s^{-23/18}  \, ds  \leq C'
\, .
\end{equation}
Finally, combining \eqr{e:echib} and \eqr{e:areabb}
gives \eqr{e:nocurv}.
\end{proof}

\begin{figure}[htbp]
    \setlength{\captionindent}{20pt}
    \begin{minipage}[t]{0.5\textwidth}
    \centering\input{pl21.pstex_t}
    \caption{A stable $\Gamma$ satisfying i)--iii): $\Gamma_0 \subset \Gamma$ is
    a disk with geodesics $\gamma_1 , \gamma_2 \subset \partial \Gamma_0$ which
    are in the middle sheets of multi-valued graphs $\Sigma_1$, $\Sigma_2$.}
    \label{f:iii}   \end{minipage}%
    \begin{minipage}[t]{0.5\textwidth}
    \centering\input{pl22.pstex_t}
    \caption{Lemma \ref{l:charc}:  A chord-arc property for a stable $\Gamma$
    satisfying i)--iii).}
    \label{f:charc}   \end{minipage}
\end{figure}

The next corollary gives upper and lower  bounds for the areas of
tubular neighborhoods in a $\Gamma$ which satisfies i)--iii)
below; see fig. \ref{f:iii}.  ($\Gamma_t (\partial)$ is the component
of $B_t \cap \Gamma$ containing $B_{t} \cap \partial \Gamma$.)\\

\noindent
 i) $\Gamma \subset B_{2R}$
is a stable embedded minimal surface, $\partial \Gamma
\subset B_{1/4} \cup
\partial B_{2R}$, $B_{1/4} \cap \partial \Gamma$ is connected,
$\Gamma_0 \subset \Gamma$ is a disk,
 $\partial
\Gamma_0 = \gamma_1 \cup \gamma_2 \cup \sigma_1 \cup \sigma_2$
($\gamma_i:[0,\ell_i] \to \Sigma$ is a geodesic),
$\gamma_i(0) \in \sigma_1 \subset \Gamma_1(\partial) $, and
$\gamma_i \perp \sigma_1$;\\
ii) $\Sigma_1 , \Sigma_2 \subset \Gamma$ are disjoint (multi-valued) graphs
over domains containing
$S_{1/2,R}^{-2\pi,2\pi}$ of $u_1 , u_2$ satisfying \eqr{e:wantit},
 $w_i < 0$, $\Sigma_i \subset \cC_{8 \, \pi \, \epsilon,
R}(u_i(1,0))$, $\partial \Sigma_i \subset \partial \cC_{8 \, \pi
\, \epsilon,R}(u_i(1,0))$, and  $\gamma_i \subset
(\Sigma_i)_{3/4,R}^{-\pi , \pi}$ ;\\
iii) $\dist_{\Gamma}
(\gamma_i (t) , \Gamma_1 (\partial)) = t$ for $0 \leq t
\leq \ell_i$, $\ell_i \geq R -1$, and
$\dist_{\Gamma} (\sigma_2 , \Gamma_1 (\partial)) \geq R-1$.

We first show  that intrinsic and
extrinsic distances to $\sigma_1$ are
roughly equivalent (see fig. \ref{f:charc}):

\begin{Lem}     \label{l:charc}
There exists $C_d>1$ so if i)--iii) hold and $R > C_d$,
 then $B_{R/C_d} \cap \sigma_2 = \emptyset$ and
$B_{t/C_d} \cap \Gamma_0 \subset \an_{t} (\sigma_1,\Gamma_0)$ for
 $C_d < t < R$.
\end{Lem}

\begin{proof}
Both of these assertions follow easily from stability together
with the assumption that $\Gamma$ contains multi-valued graphs.
Namely, suppose that either one failed.  It follows easily that
there exists a point in $\Gamma$ which is extrinsically much
closer to the origin than its intrinsic distance to the inner
boundary of $\Gamma$.  This easily implies by stability that
$\Gamma$ contains a large almost flat graph over a disk centered
at the origin which easily contradicts that $\Gamma$ contains
multi-valued graphs since these would be forced to spiral into the
almost flat graph.  We will now make this argument precise.

 Fix $C_d > 1$ to be chosen. We show first that
 $B_{t} \cap
\Gamma \subset \cT_{C_d t}(B_{1/4} \cap \partial \Gamma)$ for $1 < t < R/C_d$.
 Suppose that $y \in B_{R/C_d} \cap
\Gamma$.
Fix $C > 2$ and $\delta > 0$ to be chosen. Since $\Gamma$ is
stable, \cite{Sc}, \cite{CM2} give $C_d' = C_d' (C , \delta)$ so
that if $\dist_{\Gamma} (y, \partial \Gamma) > C_d' \, (1+|y|)$,
then
 $\cB_{C_d' \, (1+|y|)}(y)$ contains
a graph $\Gamma_y$ with gradient $\leq \delta$ over a disk $B_{C
\, (1+|y|)}(y) \cap P_y$, where $P_y \subset \RR^3$ is the plane
tangent to $\Gamma$ at $y$.  Since $\Gamma$ is embedded (and since
$\Gamma$ contains a multi-valued graph $\Sigma_1$ around
$\gamma_1$ with $\gamma_1(0) \in B_1$), we can choose $C , \delta$
so $\Gamma$ would then be forced to spiral into $\Gamma_y$. This
is impossible since $\Gamma$ is compact.  Since $\partial \Gamma \subset
\Gamma_1 (\partial) \cup \partial B_{2R}$, it follows that
$B_{t} \cap
\Gamma \subset \cT_{2C_d't}(B_{1/4} \cap \partial \Gamma)$ for $1 < t < R/C_d$.
Combining this and iii) gives
$B_{(R-1)/(2C_d')} \cap \sigma_2 = \emptyset$.

Suppose that $y \in B_{R/C_d} \cap \Gamma_0$ so (by the first part)
we get $y' \in \partial \Gamma_0$ with
\begin{equation}    \label{e:closepts}
    \dist_{\Gamma_0} (y,y') + \dist_{\Gamma} (y' , B_{1/4} \cap \partial \Gamma
    ) \leq C_d' \, (1+|y|) < R \, .
\end{equation}
In particular, $y' \in
\sigma_1 \cup \gamma_1 \cup
\gamma_2$.
 Since $\dist_{\Gamma} (\gamma_i (t) , \Gamma_1 (\partial)) = t$, we get
\begin{equation}    \label{e:cpts}
    \dist_{\sigma_1 \cup \gamma_i} (y' , \sigma_1) \leq
     C_d' \, (1+|y|) \, ,
\end{equation}
so  $\dist_{\Gamma_0} (y, \sigma_1) \leq 2 \, C_d' \, (1+ |y|)$ and
 the lemma  follows.
\end{proof}

\begin{Cor} \label{c:area}
Given $\epsilon , C_I >0$, there exists $C_4 > 0$ so if i)--iii)
hold and $R^{3/4} > 12 C_d$, then for $\e < t \leq R^{3/4}/4 - 1$
\begin{equation}  \label{e:smallbig}
    C_4 \,  \log^2 t \leq t^{-2} \Area
    (  \an_t (\sigma_1 , \Gamma_0))
     \leq  \left(1 + \int_{\an_{C_I}(\sigma_1 , \Gamma_0)} (1+ |A|^2) +
 \int_{\sigma_1}
   (1+|k_g|)
    + \log \log t \right)/C_4  \, .
\end{equation}
\end{Cor}

\begin{proof}
Since $\sigma_1 \subset \Gamma_1(\partial)$,
i) and iii) imply (A) with $C_0 =0$, (C), and (D) with $\ell = R-1$.
 Using Corollary \ref{c:bdchi} on $\Sigma_1 , \Sigma_2$,
we can define $\chi$ on $\{ x_1^2 + x_2^2 \leq R^{3/2}/4 \} \cap \Gamma$
which vanishes on $\gamma_1 , \gamma_2$, is one on
$\{ x_1^2 + x_2^2 \leq R^{3/2}/4 \} \cap \Gamma \setminus ( \Sigma_1 \cup \Sigma_2)$,
and satisfies \eqr{e:nocurv} (with double the constant).  Since
$\an_t (\sigma_1 , \Gamma_0) \subset \{ x_1^2 + x_2^2 \leq R^{3/2}/4 \}$,
inserting \eqr{e:nocurv} into Lemma
\ref{l:stable} (and scaling so $C_I \to 1$) gives
the second inequality in \eqr{e:smallbig}.

By Lemma \ref{l:cone1}, $\Sigma_1 , \Sigma_2$ each contain a
(multi-valued) graph over \eqr{e:cgd}. Suppose now that $\e < t <
R^{3/4}/(4C_d)$. By
 Lemma \ref{l:charc},
$\{1 < x_1^2 + x_2^2 \leq t^2 \} \cap \Gamma \subset B_{2t} \cap
\Gamma \subset \an_{2 C_d t}(\sigma_1 , \Gamma)$ and $B_{2t} \cap
\sigma_2 = \emptyset$ (by iii)).  Since $\sigma_1 \subset B_1$,
$\gamma_i \subset (\Sigma_i)_{3/4,R}^{-\pi , \pi}$, and  $\Sigma_1
\cap \Sigma_2 = \emptyset$, it then follows easily that $\an_{2
C_d t} (\sigma_1 , \Gamma_0)$ contains one component of $\{1 <
x_1^2 + x_2^2 \leq t^2 \} \cap \Sigma_1 \setminus
(\Sigma_1)_{1,t}^{-\pi , \pi}$.  The
 first inequality in \eqr{e:smallbig} follows immediately
(after possibly decreasing $C_4 > 0$).
\end{proof}

\begin{proof}
(of Theorem \ref{t:remsing}).  Rescale so that $r_0 =1$. Set $\hat
\Gamma=\Gamma\setminus \Gamma_{1}(\partial)$ so (since $\Gamma$ is
topologically an annulus) $\partial \hat \Gamma=\sigma\cup\hat
\sigma$ where $\sigma\subset
\partial B_{1}$, $\hat \sigma \subset \partial B_R$ are the two connected
components of $\partial \hat{\Gamma}$, and $\partial_{\nn}|x|\geq
0$ on $\sigma$ (where $\nn$ is the inward normal to $\partial \hat
\Gamma$).

By Theorem \ref{t:longstable} we need only prove that (2) does not
happen for $\hat \Gamma$. Suppose it does; we will obtain a
contradiction. The key point will be to find two oppositely
oriented multi-valued graphs in $\Gamma$ which have fixed bounded
distance between them and then apply Corollary \ref{c:area} for
$t$ sufficiently large to get a contradiction.

Fix (ordered) points $z_1 , \dots , z_m \in \sigma$ so $\sigma
\setminus \{ z_1 , \dots , z_m \}$ has components $\{ \sigma_{z_1}
, \dots , \sigma_{z_m} \}$ where $\partial \sigma_{z_i} = \{ z_i ,
z_{i+1} \}$ (set $z_{m+1} = z_1$) and ${\text{Length}}
(\sigma_{z_i}) \leq 1$. By Theorem \ref{t:longstable} (and the
discussion surrounding \eqr{e:wantit}), $\Gamma$ contains
$3$-valued graphs $\Sigma_{z_i}$ of $u_{z_i}$ satisfying
\eqr{e:wantit} over $D_{R/\omega} \setminus D_{\omega}$ (after a
rotation of $\RR^3$; a priori this rotation may depend on $z_i$)
and with $\dist_{ \hat{\Gamma} } ( z_i , (\Sigma_{z_i})_{\omega,
\omega}^{0,0}) < d_0$. Combining this with Corollary
\ref{c:longstable}, we get
 $3$-valued graphs $\{ \Sigma_{z_i} \}$,
geodesics $\gamma_{z_i} :[0,\ell_{z_i}]\to (\Sigma_{z_i})_{\omega
, R^{1/2} }^{-\pi , \pi}$ with $\gamma_{z_i} (0)\in \partial
B_{\lambda \omega}$, $\dist_{\Gamma}(\gamma_{z_i} (t),
\Gamma_{\lambda \omega} (\partial) )=t$ for $0\leq t\leq
\ell_{z_i}$,   and $\gamma_i (\ell_{z_i})\subset \Gamma \setminus
B_{R^{1/2}/3 }$. After possibly increasing $\lambda$, we can
assume that $\lambda \omega > 2 d_0 + 2$.  Hence,
 the curves in $\hat{\Gamma}$ from $z_i$ to
$(\Sigma_{z_i})_{\omega, \omega}^{0,0}$ given by Theorem
\ref{t:longstable} are contained in $B_{\lambda \omega/2}$.
Therefore, since $(\Sigma_{z_i})_{\omega, \lambda
\omega}^{-3\pi,3\pi}$ is a graph,   we can choose curves
$\eta_{z_i} \subset \Gamma_{\lambda \omega}(\partial)$ from
$\gamma_{z_i}(0)$ to $z_i$ with length $\leq 2 \lambda \omega +
4\pi \omega$ and so $\eta_{z_i}\setminus B_{\lambda \omega/2}$ is
simple with $\int_{\eta_{z_i}\setminus B_{\lambda \omega/2}} |k_g|
\leq C$.

It follows immediately from embeddedness that the $\Sigma_{z_i}$'s
are graphs over a common plane.   From the gradient estimate
(which applies because of estimates for stable surfaces of
\cite{Sc}, \cite{CM2}),   each component of  $\Gamma$ intersected
with a concave cone is also a multi-valued graph.  Since
 $\partial B_{\lambda \omega} \cap
 \partial \Gamma_{\lambda \omega}(\partial)$ is a closed curve,
 it must pass between the sheets of each $\Sigma_{z_i}$.
It is now easy to see that each $\Sigma_{z_i}$ contains an
oppositely oriented multi-valued graph $\hat{\Sigma}_{z_i}$
between its sheets (i.e., $\nn_{\Gamma}$ points in almost opposite
directions on $\Sigma_{z_i}$ and $\hat{\Sigma}_{z_i}$).
Furthermore, since Lemma \ref{l:charc} bounds the distance in
$\hat{\Gamma}$ from $\hat{\Sigma}_{z_i}$ to $\sigma$, we can
assume that two of the $\Sigma_{z_i}$'s are oppositely oriented.
 We can therefore choose two consecutive $3$-valued graphs,
$\Sigma_{z_j}$, $\Sigma_{z_{j+1}}$, which are  oppositely
oriented; rename these $\Sigma_1$, $\Sigma_2$ (and similarly the
corresponding $\gamma_1$, $\gamma_2$, $\ell_1$, $\ell_2$).

By replacing $B_{\lambda \omega/2} \cap (\sigma_{z_j} \cup \eta_{z_j} \cup
\eta_{z_{j+1}})$  with a broken
geodesic and finding a simple subcurve as in Lemma \ref{l:broken},
we get a simple curve $\sigma_1 \subset  \Gamma_{\lambda
\omega}(\partial) \setminus \Gamma_{7/8}(\partial)$
from $\gamma_1(0)$ to $\gamma_2(0)$ with
\begin{equation}    \label{e:getnow0}
   \int_{\sigma_1}
   (1+|k_g|) \leq C_a \, .
\end{equation}
Furthermore, since $\sigma_1 \subset \Gamma_{\lambda
\omega}(\partial)$),
  $\dist_{\Gamma} (\gamma_i(t) , \sigma_1) = t$
for $0\leq t \leq \ell_i$.
 Let $\Gamma_0$ be the
component of $\Gamma_{R^{1/2}/3}(\partial) \setminus ( \sigma_1 \cup \gamma_1
\cup \gamma_2)$ which does not contain $\Gamma_{7/8}(\partial)$;
set $\sigma_2 = \partial \Gamma_0 \setminus ( \sigma_1 \cup \gamma_1
\cup \gamma_2)$.
  It follows
that $\Gamma_0$ is a disk and $\dist_{\Gamma} (\Gamma_0 , \partial
\Gamma) \geq 5/8$. Since $(\Sigma_{z_i})_{\omega, \lambda
\omega}^{-3\pi,3\pi}$ is a graph, we can perturb $\sigma_1$ near
$\gamma_1 (0) , \gamma_2(0)$ to arrange that $\sigma_1 \perp
\gamma_1$ and $\sigma_1 \perp \gamma_2$ and so $\sigma_1$ still
satisfies \eqr{e:getnow0} with a slightly larger constant $C_a$.
 Combining \eqr{e:getnow0} and estimates
for stable surfaces of \cite{Sc}, \cite{CM2}, we get
\begin{equation}    \label{e:getnow}
  \int_{ \an_{1/8}({\sigma}_1 , {\Gamma}_0)} (1+|A|^2) +
 \int_{{\sigma}_1}
   (1+|k_g|) \leq C_b \, .
\end{equation}
Hence, (after rescaling) ${\Gamma}_0 , \Gamma , \Sigma_1 ,
\Sigma_2 , {\gamma}_1 , {\gamma}_2 , {\sigma}_1$ satisfy i)--iii).
To get ii), we use \cite{Sc}, \cite{CM2} and the gradient estimate
to extend $\Sigma_1 , \Sigma_2$ as multi-valued graphs inside the
cones $\cC_{8 \, \pi \, \epsilon, R^{1/2}/4}(u_i(1,0))$; the
opposite orientation guarantees that $\Sigma_1 \cap \Sigma_2 =
\emptyset$. Corollary \ref{c:area} and \eqr{e:getnow} give for
$C_5 < t < R^{3/8}/C_5$
\begin{equation}  \label{e:smallbig2}
    C_4 \,  \log^2 t \leq t^{-2} \Area
    (  \an_t ({\sigma}_1 , {\Gamma}_0))
     \leq  \left(1 + C_b
    + \log \log t \right)/C_4  \, .
\end{equation}
This gives
 the desired
contradiction for $R$  large, completing the proof.
\end{proof}

\part{Nearby points with large curvature}       \label{p:p3}

In this part, we extend Theorem \ref{t:remsing} (proven for stable
surfaces) to surfaces with {\it extrinsic} quadratic curvature
decay $|A|^2 \leq C \, |x|^{-2}$.  As mentioned in the introduction, this
extension is needed in both \cite{CM5} and \cite{CM6}.  In \cite{CM5} it
is used for disks to get points of large curvature nearby and on each side
of a given point with large curvature (in particular it is used to show
that such points are not extrinsically isolated).

Stability was used in the proof of Theorem \ref{t:remsing} for two
purposes: (a) To show {\it intrinsic} quadratic curvature decay.
(b) To bound the total curvature using the stability inequality.
To get the extension to the extrinsic quadratic curvature decay
case, we will deal with (a) and (b) separately in the next two
sections. To get (a), we relate extrinsic and intrinsic distances
(i.e., we show a ``chord-arc'' property). For (b), we follow
section $2$ of \cite{CM4} to decompose a surface with quadatric
curvature decay into disjoint almost stable subdomains and a
``remainder'' with quadratic area growth.

  For applications of the results of this part in
 \cite{CM5}, $\Sigma$ will be a disk and hence $\partial \Sigma_{0,t}$
is connected for all $t$  (here, and elsewhere,  if $0\in \Sigma$,
then $\Sigma_{0,t}$ denotes the component of $B_t \cap \Sigma$
containing $0$).  However, in \cite{CM6}, when we apply the
results here to deal with the first possibility in  ``4).'' of
Theorem \ref{t:uniplanar} (i.e., the analog of the genus one
helicoid), $\Sigma$ is no longer a disk but $\partial \Sigma$ is
still connected (which is assumed in many of the results below).

\section{Relating intrinsic and extrinsic distances}

In this section,  $0 \in \Sigma \subset B_{R}$ is an embedded
minimal surface, $\partial  \Sigma \subset
\partial B_{R}$,  $|A|^2 \leq C_1^2 |x|^{-2}$ on $\Sigma
\setminus B_1$, and $\partial \Sigma_{0,t}$ is connected for $1
\leq t \leq R$.

 The next lemma shows that  only one component of $B_{C_b}
\cap \Sigma$ intersects $B_2$. The second lemma bounds the radius
of the intrinsic  tubular neighborhood of $B_2 \cap \Sigma$
containing this component.  Combining these iteratively (on
decreasing scales) in Corollary \ref{l:extcc} gives the
``chord-arc'' property needed to establish (a).

\begin{Lem} \label{l:extca}
Given $C_1$, there exists $C_b$ so
 if $\Sigma_{0,1}$ is not a graph, then $B_2 \cap \Sigma \subset
 \Sigma_{0,C_b}$.
\end{Lem}

\begin{proof}
Suppose that $\Sigma_1 , \Sigma_2$ are disjoint components of
$B_{C_b} \cap \Sigma$ with $B_2 \cap \Sigma_i \ne \emptyset$.
It follows that there is a component $\Omega$ of
$B_{C_b} \setminus \Sigma$ and a segment $\eta \subset B_2 \setminus
\Sigma$ so that $\partial \Sigma_{0,C_b}$ is linked with $\eta$
in $\Omega$ (cf. lemma $2.1$ in \cite{CM9}).  Since $\Omega$ is
mean convex, we can solve the Plateau problem as in \cite{MeYa2}
to get a stable minimal surface $\Gamma
\subset \Omega$ with
$\partial \Gamma = \partial \Sigma_{0,C_b}$.
The linking implies that $B_2 \cap \Gamma \ne
\emptyset$.
The curvature estimates of
\cite{Sc}, \cite{CM2} then give
a graph $\Gamma_0 \subset
\Gamma$ of a function $u_0$ over $D_{C_b/C}$ (after a
rotation) with $|u_0(z)| \leq |z|$. By corollary $1.14$ of
\cite{CM8} (applied with $w=0$),
we can assume that on $D_{C_b^{1/2}/C}$
\begin{equation}    \label{e:fromcm8}
    |\nabla u_0| (z) \leq C' \, |z|^{-5/12} \, .
\end{equation}
In particular, $\Gamma_0$ is close to a horizontal plane. The
lemma now follows from an argument used in \cite{CM9} (see also
\cite{CM10}) which we now outline: $\Sigma$ intersects a narrow
cone about $\Gamma_0$, then contains a long chain of graphical
balls (by the gradient estimate), and must then either spiral
indefinitely or close up as a graph. Namely, for $t< C_b^{1/2}/C$,
$\Sigma_{0,t}$ sits on one side of $\Gamma_0$. However, by lemma
$2.4$ of \cite{CM9} (for $t > C'$),
$\partial \Sigma_{0,t}$ contains a ``low
point'' $y_0$ (i.e., $|x_3 (y_0)| \leq \delta \, t$ with $\delta >
0$ small). The gradient estimate (since $|A|^2 \leq C_1^2
|x|^{-2}$ on $\Sigma \setminus B_1$) gives a long chain of balls
$\cB_{c\,t}(y_i)$ with $y_i \in
\partial \Sigma_{0,t} \cap \{ |x_3| \leq C' \, \delta \, t \}$ and
which is a (possibly multi-valued) graph. Since $\partial
\Sigma_{0,t}$ cannot spiral forever, this graph closes up.
By Rado's theorem (note that
no assumption on the topology is needed for this application of Rado's theorem;
cf. the proof of theorem $1.22$ in
\cite{CM4}), $\Sigma_{0,t}$ is itself a graph,
giving the lemma.
\end{proof}

\begin{Lem} \label{l:extcb}
Given $C_1 , C_b$, there exists $C_{c}$ so if
$R > C_c$, then for all $y \in \Sigma_{0,C_b}$
\begin{equation}    \label{e:extcb}
    \dist_{\Sigma} (y , B_1 \cap \Sigma) \leq C_{c} \, .
\end{equation}
\end{Lem}

\begin{proof}
Let $\tilde{\Sigma}$ be the universal cover of $\Sigma$ and $\tpi
: \tilde{\Sigma} \to \Sigma$ the covering map. With the definition
of $\delta$-stable as in section $2$ of \cite{CM4}, the argument
of \cite{CM2} (i.e., curvature estimates for $1/2$-stable
surfaces) gives
 $C > 10$ so if
$\cB_{C C_b/2 }(\tilde{z}) \subset \tilde{\Sigma}$ is $1/2$-stable
and  $\tpi(\tilde{z}) = z$, then $\tpi : \cB_{5C_b}(\tilde{z}) \to
\cB_{5C_b}(z)$ is one-to-one and $\cB_{5C_b}(z)$ is a graph with
$B_{4C_b}(z) \cap
\partial \cB_{5C_b}(z) = \emptyset$. Corollary 2.13 in \cite{CM4}
gives $\epsilon  = \epsilon (C , C_1 , C_b) > 0$ so that if $
|z_{1} - z_{2}| < \epsilon$ and $|A|^2 \leq C_1^2$ on (the
disjoint balls) $\cB_{ C C_b }(z_i)$, then each $\cB_{C C_b
/2}(\tilde{z}_{i}) \subset \tilde{\Sigma}$ is $1/2$-stable where
$\tpi (\tilde{z}_{i}) = z_i$.

We claim that there exists $n$ so $B_1 \cap \cB_{(2n+1) \, C C_b
}(y) \ne \emptyset$. Suppose not; we get a curve $\sigma \subset
\Sigma_{0,C_b}  \setminus \cT_{C C_b }(B_{1} \cap \Sigma)$ from
$y$ to $\partial \cB_{2n\,C C_b }(y)$. For $i=1, \dots , n$, fix
points $z_i \in \partial \cB_{2i\,C C_b }(y) \cap \sigma$. The
intrinsic balls $\cB_{C C_b} (z_i) \subset \Sigma \setminus B_1$
are disjoint, have centers in $B_{C_b} \subset \RR^3$, and
$\sup_{\cB_{C C_b }(z_i)} |A|^2 \leq C_1^2$. Hence, there exist
$i_1 , i_2$ with $0< | z_{i_1} - z_{i_2}| < C' \, C_b \, n^{-1/3}
< \epsilon$, and, by corollary 2.13 in \cite{CM4}, each $\cB_{C
C_b /2}(\tilde{z}_{i_j}) \subset \tilde{\Sigma}$ is $1/2$-stable
where $\tpi ( \tilde{z}_{i_j}) = {z}_{i_j}$. By \cite{CM2}, each
$\cB_{5C_b}(z_{i_j})$ is a graph with $B_{4C_b}(z_{i_j}) \cap
\partial \cB_{5C_b}({z}_{i_j}) = \emptyset$.  In particular,
$B_{C_b} \cap
\partial \cB_{5C_b}(z_{i_j}) = \emptyset$.
This contradicts that $\sigma \subset B_{C_b}$ connects
$z_{i_j}$ to $\partial \cB_{C C_b} (z_{i_j})$.
\end{proof}

\begin{Cor} \label{l:extcc}
Given $C_1$, there exists $C_{c}$ so if $\Sigma_{0,1}$ is not a
graph, $y \in B_{R/C_c} \cap \Sigma$, then
\begin{equation}    \label{e:extcc}
    \dist_{\Sigma} (y , B_1 \cap \Sigma) \leq 2 \, C_{c} \, |y| \, .
\end{equation}
\end{Cor}

\begin{proof}
Suppose $y \in B_{2^n} \setminus B_{2^{n-1}}$.  By Lemma
\ref{l:extca}, $y \in \Sigma_{0,C_b 2^{n-1}}$ where $C_b = C_b
(C_1)$.  Set $y_n = y$. Lemma \ref{l:extcb}  gives $y_{n-1} \in
B_{2^{n-1}} \cap \Sigma$ with $\dist_{\Sigma} (y_n ,y_{n-1}) \leq
C_c \, 2^{n-1}$.

We can now repeat the argument.  Namely, by Lemma \ref{l:extca},
$y_{n-1} \in \Sigma_{0,C_b 2^{n-2}}$ and then Lemma \ref{l:extcb}
gives $y_{n-2} \in B_{2^{n-2}} \cap \Sigma$ with $\dist_{\Sigma}
(y_{n-1} ,y_{n-2}) \leq C_c \, 2^{n-2}$, etc. After $n$ steps, we
get $y_0 \in B_1 \cap \Sigma$ with
\begin{equation}    \label{e:finally}
    \dist_{\Sigma} (y,y_{0}) \leq \sum_{i=1}^{n}
    \dist_{\Sigma} (y_i,y_{i-1}) \leq \sum_{i=1}^{n} C_c \, 2^{i-1}
    \leq 2 \, C_c \, |y| \, .
\end{equation}
\end{proof}

\section{A decomposition from \cite{CM4}}       \label{s:decm4}

In lemma $2.15$ of \cite{CM4}, we decomposed an embedded
minimal surface in a ball with
bounded curvature into disjoint almost stable subdomains and a
remainder with bounded area.  The same argument gives:

\begin{Lem} \label{l:deltstb}
Given $C_1$, there exists $C_d$ so: If $\Sigma \subset B_{2R}$ is
an embedded minimal surface with $\partial \Sigma \subset
\partial B_{2R} \cup B_{1/2}$, and $|A|^2 \leq C_1^2 \, |x|^{-2}$,
then there exist disjoint $1/2$-stable subdomains $\Omega_j
\subset \Sigma$ and a function $0 \leq \psi \leq 1$ on
$\Sigma$ which vanishes on
$(B_R \setminus B_1) \cap \Sigma \setminus (\cup_j \Omega_j)$ so
\begin{align} \label{e:nabd}
\Area ( \{ x \in (B_R \setminus B_1) \cap \Sigma \, | \, \psi (x)
< 1 \} ) & \leq C_d \, R^2 \, , \\ \int_{B_R \cap \Sigma} |\nabla \psi|^2 &
\leq C_d \, \log R \, . \label{e:nchi}
\end{align}
\end{Lem}

In the proof of Theorem \ref{t:remsingG} in the next section,
Lemma \ref{l:deltstb}
will be used to extend the area bounds for
stable surfaces proven in Sections \ref{s:lst} and
\ref{s:areagrwth}
(specifically those in Lemma \ref{l:stable}, Proposition \ref{c:bbl},
and Corollary \ref{c:area})
to minimal surfaces with $|A|^2 \leq C_1^2 \, |x|^{-2}$.  This is
very similar to how lemma $2.15$ of \cite{CM4} was used in
lemma $3.1$ of \cite{CM4}.

 By Lemma \ref{l:deltstb}, $\int_{B_R \cap \Sigma} |\nabla \psi|^2 +
 \int_{B_R \cap \{ \psi < 1 \} } |A|^2$  grows (in $R$) at most like
$\log R$.  We use this below in the $1/2$-stability inequality to
get the total curvature bound needed for (b).  This is used in the
proof of Theorem \ref{t:remsingG1}.

\section{Theorem \ref{t:remsingG} and a generalization}

As already mentioned, stability was used in the
proof of Theorem \ref{t:remsing} to establish (a) and (b) in the
introduction to Part \ref{p:p3}; these were extended in the two
previous sections to surfaces with a quadratic curvature bound. In
\cite{CM5} we will need the contrapositive of Theorem \ref{t:remsingG},
i.e., we will need to find points where the quadratic bound fails.
In fact, what we will really need is to find points on ``each
side'' of a multi-valued graph where this fails; this is the
following theorem:

(Here
 $u_1(r_0,2\pi) < u_2 (r_0,0) < u_1 (r_0,0)$ just says that the two graphs
 spiral together, one inside the other; cf. theorem $0.6$ in \cite{CM8}.)

\begin{figure}[htbp]
    \setlength{\captionindent}{20pt}
    \begin{minipage}[t]{0.5\textwidth}
    \centering\input{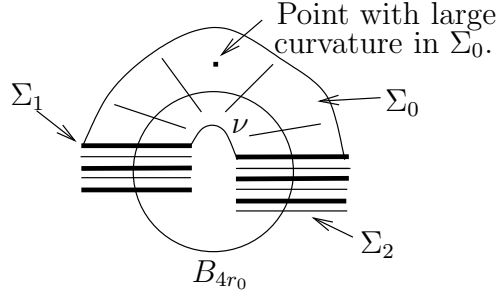}
    \caption{Theorem \ref{t:remsingG1} and Corollary \ref{c:remsingG1}
    - existence of nearby points with large curvature.}
    \label{f:nearby}   \end{minipage}
\end{figure}

\begin{Thm} \label{t:remsingG1}
See fig. \ref{f:nearby}.
Given $C_1$, there exists $C_2$ so: Let  $0 \in\Sigma \subset B_{2
C_2 \, r_0}$ be an embedded minimal surface with connected
$\partial  \Sigma \subset \partial B_{2 C_2 \, r_0}$ and $\Genus
(\Sigma_{0,r_0}) = \Genus (\Sigma)$. Suppose $ \Sigma_1 , \Sigma_2
\subset  \Sigma\cap \{x_3^2\leq (x_1^2+x_2^2)\}$ are
(multi-valued) graphs of $u_i$ satisfying \eqr{e:wantit} on
$S_{r_0, C_2 r_0}^{-2\pi,2\pi}$, $u_1(r_0,2\pi) < u_2 (r_0,0) <
u_1 (r_0,0)$, and $\nu \subset \partial \Sigma_{0,2r_0}$ a curve
from $\Sigma_1$ to $\Sigma_2$. If  $\Sigma_0$ is the component of
$\Sigma_{0,C_2 r_0} \setminus (\Sigma_1 \cup \Sigma_2 \cup \nu)$
which does not contain $\Sigma_{0,r_0}$, then
\begin{equation}    \label{e:nextone1}
    \sup_{x\in \Sigma_{0} \setminus B_{4r_0}} |x|^2 \, |A|^2 (x)
\geq 4\,C_1^2 \, .
\end{equation}
\end{Thm}

\begin{proof}
Suppose that
 \eqr{e:nextone1} fails for some $C_1$; as in the proof of Theorem
\ref{t:remsing}, we will show contradictory upper and lower bounds for
the area growth for
$C_2$ sufficiently large.

Note that  for $r_0 \leq s \leq 2C_2 r_0$,
  it follows from the
maximum principle (since $\Sigma$ is minimal) and Corollary
\ref{c:gen2} that
$\partial \Sigma_{0,s}$ is connected and $\Sigma \setminus
\Sigma_{0,s}$ is an annulus.

Note also that the gradient estimate (which applies because of the
curvature bound)  allows us to extend each $\Sigma_i$ (inside
$\Sigma_0$) as a graph of $u_i$ over $\partial D_{\rho}$  as long
as $|u_i (\rho,\theta) - u_i (\rho , [\theta])| \leq C_g \, \rho$,
where $\theta - [\theta] \in 2\pi \ZZ$ and $0\leq [\theta] \leq
2\pi$.  By corollary 1.14 of \cite{CM8},   the curvature of
$\Sigma_i$ decays faster than quadratically. Combining these (and
increasing the inner radius), we can assume that each $\Sigma_i$
extends (inside $\Sigma_0$) as a graph until it leaves a cone $\{
x_3^2 \leq \Lambda^2 (x_1^2 + x_2^2) \}$ for some small $\Lambda >
0$. Moreover, these extended multi-valued graphs must stay
disjoint since $u_1(r_0,2\pi) < u_2 (r_0,0) < u_1 (r_0,0)$.

We next choose the inner boundary curve where we argue as in
Theorem \ref{t:remsing}.  By Lemma \ref{l:extca}, $B_{4r_0} \cap
\Sigma \subset\Sigma_{0,2C_b r_0}$. In particular, $\partial
\Sigma_{0,2C_b r_0}$ separates $B_{4r_0} \cap \Sigma$ from
$\partial \Sigma$.  We can therefore replace $\nu$ with a segment
of $\partial \Sigma_{0,2C_b r_0}$ from $\Sigma_1$ to $\Sigma_2$ so
(for the new $\Sigma_0$)
\begin{equation}    \label{e:gettoassume}
    \sup_{x\in \Sigma_{0}} |x|^2 \, |A|^2 (x)
\leq 4\, \bar{C}_1^2 \, .
\end{equation}
By Corollary \ref{l:extcc} (the ``chord-arc'' property), intrinsic
and extrinsic distances to $B_{4r_0} \cap \Sigma$ are compatible.
Hence, we get
\begin{equation}    \label{e:nof1}
    \sup_{x\in \Sigma_{0}} \, \dist_{\Sigma}^{2} (x, B_{4r_0} \cap \Sigma ) \,
        |A|^2 (x) \leq C_3  \, .
\end{equation}

The proof of Theorem \ref{t:remsing} now applies with two changes
(and the  minor modifications which result):\\
(a') The curvature estimates for stable surfaces of \cite{Sc},
\cite{CM2} are replaced with
\eqr{e:nof1}.\\
(b') The total curvature bound from the stability
inequality in \eqr{e:intcut}
 is replaced with
the bound using Lemma \ref{l:deltstb} and the $1/2$-stability
inequality (cf.  lemma $3.1$ of
\cite{CM4}).

Namely, using (a') and (b'),
the proof of Theorem \ref{t:longstable} extends from stable surfaces
to surfaces satisfying \eqr{e:nof1}
(with (b') being used in Lemma \ref{l:stable} and
Proposition \ref{c:bbl}  exactly as in \cite{CM4}).
 It follows that each $z$ in (the new) $\nu$ is a fixed bounded
distance from a multi-valued graph (either $\Sigma_1 , \Sigma_2$
or a new multi-valued graph in between).  Hence, as in the proof
of Theorem \ref{t:remsing}, we can choose two consecutive
multi-valued graphs which are  oppositely oriented; let $\sigma_1$
be the curve connecting these. Next, (b')  contributes a new $C_4
\, t^2 \, \log t$ term to the upper bound for the  area of a
sector $\an_t (\sigma_1)$ in the upper bound for the area in
Corollary \ref{c:area} where $C_4$ does not depend on $\sigma_1$
(see the last paragraph of Section \ref{s:decm4}). However, since
the lower bound for the area is on the order of $t^2 \, \log^2 t$,
we get the desired contradiction as before.
\end{proof}

In \cite{CM5}, we will use the special case of
Theorem \ref{t:remsingG1} where $\Sigma$ is a disk:

\begin{Cor} \label{c:remsingG1}
See fig. \ref{f:nearby}.
Given $C_1$, there exists $C_2$
so: Let  $0 \in\Sigma \subset B_{2 C_2 \, r_0}$ be an embedded
minimal disk. Suppose $ \Sigma_1 , \Sigma_2 \subset  \Sigma\cap
\{x_3^2\leq (x_1^2+x_2^2)\}$ are graphs of $u_i$ satisfying
\eqr{e:wantit} on $S_{r_0, C_2 r_0}^{-2\pi,2\pi}$, $u_1(r_0,2\pi)
< u_2 (r_0,0) < u_1 (r_0,0)$, and $\nu \subset \partial
\Sigma_{0,2r_0}$ a curve from $\Sigma_1$ to $\Sigma_2$. Let
$\Sigma_0$ be the component of $\Sigma_{0,C_2 r_0} \setminus
(\Sigma_1 \cup \Sigma_2 \cup \nu)$ which does not contain
$\Sigma_{0,r_0}$.  Suppose  either
 $\partial \Sigma \subset \partial B_{2 C_2 \, r_0}$ or
$\Sigma$ is stable and
$\Sigma_0$  does not intersect $\partial \Sigma$. Then
\begin{equation}    \label{e:nextone1c}
    \sup_{x\in \Sigma_{0} \setminus B_{4r_0}} |x|^2 \, |A|^2 (x)
\geq 4\,C_1^2 \, .
\end{equation}
\end{Cor}

\begin{proof}
Since $\Sigma$ is a disk,
$\partial \Sigma$ is connected and
$\Genus (\Sigma_{0,r_0}) =\Genus(\Sigma) = 0$. Hence,
Theorem \ref{t:remsingG1}
gives the corollary when $\partial \Sigma \subset \partial B_{2 C_2 \, r_0}$.

When $\Sigma$ is stable and $\Sigma_0$  does not intersect
$\partial \Sigma$, then $\Sigma_1 , \Sigma_2$ each extend inside
cones in at least one direction as multi-valued graphs.  This
gives essentially half of the multi-valued graphs $\Sigma_1
,\Sigma_2$ used in Section \ref{s:areagrwth} which is all that is
needed in the proof of Theorem \ref{t:remsing}. The corollary now
follows easily from the proof of Theorem \ref{t:remsing} (with
$\Sigma_1 , \Sigma_2$ causing the same modifications as in Theorem
\ref{t:remsingG1}).
\end{proof}

Note that if $C_1$ is large, then \eqr{e:nextone1c} would
contradict the curvature estimate for stable
surfaces of \cite{Sc}, \cite{CM2}.  In \cite{CM5},
we will apply Corollary \ref{c:remsingG1} in this way, showing that
 such a stable $\Sigma$ does not exist.

In \cite{CM5}, we will also use the other case of
Corollary  \ref{c:remsingG1}, where $\Sigma$ is not assumed to be
stable,  to get points
of large curvature ``metrically'' on each side of the multi-valued graph
$\Sigma_1$.
Namely, note first that the curve $\partial \Sigma_{0,2r_0} \setminus \nu$ in
 Corollary \ref{c:remsingG1}  has
the same properties as $\nu$.  In \cite{CM5}, $\nu$
(and hence also $\Sigma_0$) will be
on one side of $\Sigma_1 ,\Sigma_2$ while
$\partial \Sigma_{0,2r_0} \setminus \nu$ is on the other.
Applying
Corollary \ref{c:remsingG1} to each of these will give points of
large curvature ``topologically '' on each side of $\Sigma_1 ,\Sigma_2$.

In fact, we will see in \cite{CM5} that if an embedded minimal disk
$\Sigma$ contains one multi-valued graph $\Sigma_1$, then it will
 contain
a second multi-valued graph  $\Sigma_2$ which spirals together with $\Sigma_1$
(``the other half''). We will also see there that
$\partial \Sigma_{0,Cr_0} \setminus (\Sigma_1 \cup \Sigma_2)$
has exactly two components $\nu_{\pm}$;   it follows easily
 that we can assume $\nu_+$ is above and $\nu_-$ is
below $\Sigma_1$.   Applying
Corollary \ref{c:remsingG1} to both $\nu_{\pm}$ will give points of
large curvature ``metrically'' on each side of $\Sigma_1$.

\begin{proof}
(of Theorem \ref{t:remsingG}). It suffices to show that if $\Area
(\Sigma_{0,r_0} ) >  C_3 \, r_0^2$, then
 \eqr{e:atcbG} fails.

Note that  for $r_0 \leq s \leq R$,
  it follows from the
maximum principle (since $\Sigma$ is minimal) and Corollary
\ref{c:gen2} that $\partial \Sigma_{0,s}$ is connected and $\Sigma
\setminus \Sigma_{0,s}$ is an annulus.

The proof is now virtually identical to the proof of  Theorem
\ref{t:remsingG1} except that it simplifies since we no longer
keep track of the two sides and (1) in (analog of) Theorem
\ref{t:longstable} becomes $\Area (\Sigma_{0,r_0} ) \leq  C_3' \,
r_0^2$.
\end{proof}

\end{document}